\input amstex
\documentstyle{amsppt}
\document
\magnification=1200
\NoBlackBoxes
\nologo
\vsize16cm



\centerline{\bf ITERATED INTEGRALS OF MODULAR FORMS}

\medskip

\centerline{\bf AND NONCOMMUTATIVE MODULAR SYMBOLS}

\medskip

\centerline{\bf Yuri I. Manin}

\medskip

\centerline{\it Max--Planck--Institut f\"ur Mathematik, Bonn, Germany,}

\centerline{\it and Northwestern University, Evanston, USA}

\bigskip

{\bf Abstract.} The main goal of this paper is to study
properties of the iterated integrals of modular forms
in the upper halfplane, eventually multiplied by $z^{s-1}$,
along geodesics connecting two cusps. This setting
generalizes simultaneously the theory of modular symbols
and that of multiple zeta values.

\bigskip 


\bigskip

\centerline{\bf \S 0. Introduction and summary}

\medskip

This paper was inspired by two sources: theory of multiple zeta
values on the one hand (see [Za2]), and theory of modular symbols
and periods of cusp forms, on the other ([Ma1], [Ma2], [Sh1]--[Sh3],
[Me]).
Roughly speaking, it extends the theory of  periods of modular forms
replacing integration along geodesics in the upper complex
half--plane by iterated integration. Here are some details.

\medskip

{\bf 0.1. Multiple zeta values.} They are the numbers given
by the $k$--multiple Dirichlet series
$$
\zeta (m_1,\dots ,m_k)=\sum_{0<n_1<\dots <n_k}
\frac{1}{n_1^{m_1}\dots n_k^{m_k}}
\eqno(0.1)
$$
which converge for all integer $m_i\ge 1$ and $m_k>1$,
or equivalently by the $m$--multiple iterated integrals,
$m=m_1+\dots +m_k$,
$$
\zeta (m_1,\dots ,m_k) = \int_0^1 \frac{dz_1}{z_1}
\int_0^{z_1} \frac{dz_2}{z_{2}}\int_0^{z_2} \dots
\int_0^{z_{m_{k-1}}} \frac{dz_{m_k}}{1-z_{m_k}}\dots
\eqno(0.2) 
$$
where the sequence of differential forms in the
iterated integral consists of consecutive subsequences of the
form $\dfrac{dz}{z},\dots , \dfrac{dz}{z},\dfrac{dz}{1-z}$
of lengths $m_k, m_{k-1}, \dots , m_1.$

\smallskip

Easy combinatorial considerations allow one to express
in two different ways products 
$ \zeta (l_1,\dots ,l_j)\cdot \zeta (m_1,\dots ,m_k)$
as linear combinations of multiple zeta values.

\smallskip

If one uses for this the integral representation (0.2),
one gets a sum over shuffles which enumerate the 
simplices of highest dimension occurring in the
natural simplicial decomposition of the product
of two integration simplices.

\smallskip 

If one uses instead (0.1), one gets sums
over shuffles with repetitions which enumerate some
simplices of lower dimension as well. 

\smallskip

These relations and
their consequences are called double shuffle relations.
Both types of relations can be succinctly written down in terms
of formal series on free noncommutating generators.
One can include in these relations regularized multiple zeta values
for arguments where the convergence of (0.1), (0.2) fails.

\smallskip

For a very clear and systematic exposition of these results,
see [De] and [Ra1], [Ra2].

\smallskip

In fact, the formal generating series for (regularized) iterated integrals
(0.2) appeared in the famous Drinfeld paper [Dr2], 
essentially as {\it the Drinfeld associator}, and more relations
for multiple zeta values were implicitly deduced there.
The question about interdependence of (double) shuffle 
and associator relations does not seem to be settled at the moment of writing this: cf. [Ra3]. The problem of completeness of these systems of relations
is equivalent to some difficult transcendence questions.

\smallskip

Multiple zeta values are interesting, because they and their
generalizations appear
in many different contexts involving mixed Tate motives ([DeGo], [T]),
deformation quantization ([Kon]), knot invariants etc.

\medskip

{\bf 0.2. Modular symbols and periods of modular forms.} Let
$\Gamma$ be a congruence subgroup of the modular group
acting upon the union $\overline{H}$ of the upper complex 
half--plane $H$ and the set of
cusps $\bold{P}^1(\bold{Q})$. 

\smallskip

The quotient $\Gamma \setminus \overline{H}$ is the modular curve $X_{\Gamma}$.
Differentials of the first kind on $X_{\Gamma}$ lift to the
cusp forms of weight 2 on $H$ (multiplied by $dz$).

\smallskip

The modular symbols $\{\alpha ,\beta\}_{\Gamma}\in H_1(X_{\Gamma},\bold{Q})$,
where $\alpha ,\beta \in \bold{P}^1(\bold{Q})$,
were introduced in [Ma1] as linear functionals on the
space of differentials of the first kind obtained by lifting and
integrating. The fact that one lands in $H_1(X_{\Gamma},\bold{Q})$
and not just $H_1(X_{\Gamma},\bold{R})$ is not obvious. It was proved in [Dr1] 
by refining a weaker argument given in [Ma1]. This is equivalent to the statement
that difference of any two cusps in $\Gamma$ has finite order
in the Jacobian, or else that the mixed Hodge structure on
$H^1(X_{\Gamma}\setminus \{cusps\},\bold{Q})$ is split (cf. [El]).
One of the basic new insights of [Ma1] consisted in the realization that
studying the action of Hecke operators on modular symbols
one gets new arithmetic facts about periods and Fourier
coefficients of cusp forms of weight two.

\smallskip

The further generalizations of modular symbols proceeded,
in particular, in the following
directions.

\smallskip

(a) In [Ma2] it was demonstrated that the same technique
applies to the integrals of cusp forms of higher weight,
eventually multiplied by polynomials in $z$,
producing the similar information about their periods and Fourier
coefficients. In principle, such integrals cannot be pushed down
to $X_{\Gamma}$, but they can be pushed down to 
the appropriate Kuga--Sato varieties over $X_{\Gamma}$
that is, relative Cartesian powers of the universal elliptic curve.
In this way, modular symbols of higher weight can be interpreted as
rational homology classes of middle dimension
of  Kuga--Sato varieties: cf. [Sh1]--[Sh3].

\smallskip

(b) Pushing down an oriented geodesic connecting two cusps
in $\overline{H}$ to $X_{\Gamma}$, we get a singular
chain with a boundary in cusps, which is a relative cycle
modulo cusps with integral coefficients. This is the viewpoint of [Me].
Hence it is more natural to consider
the relative/non--compact version of modular symbols,
and allow integration of the Eisenstein series, that is, differential
forms of the third kind with poles at cusps as well.
The same remark applies to the modular symbols of higher weight.

\smallskip

This refinement appears as well in the study of the
``noncommutative boundary'' of the modular space,
that is, the (tower of) space(s) $\Gamma\setminus \bold{P}^1(\bold{R})$,
cf. [MaMar]. Namely, it turns out that the relative 1--homology
modulo cusps (and additional groups of similar nature)
can be interpreted as (sub)groups of the $K$--theory
of the noncommutative boundary.

\smallskip

In this paper I suggest a generalization in the third direction, namely

\smallskip

(c) The study of iterated integrals of cusp forms and
Eisenstein series, eventually multiplied 
by a power of $z$, along geodesics connecting two cusps.
Some of these integrals can be pushed down to $X_{\Gamma}$
and thus produce a de Rham version of modular symbols which assigns
iterated (eventually regularized) periods to the elements of
the fundamental groupoid of $(X_{\Gamma}, \{cusps\})$
instead of its 1-homology group.
One may call them {\it noncommutative modular symbols.}

\smallskip

Other integrals can only be pushed down to the Kuga--Sato
varieties, or preferably, to some (covers of the) moduli
spaces $\overline{M}_{1,n}$, in the same vein as it was done for
multiple zeta values and $\overline{M}_{0,n}$ in [GoMa].
The related geometry deserves further study, both for integrands
related to cusp forms and to Eisenstein series. 

\smallskip

Notice in conclusion that the discussion above implicitly
referred only to the case
of $SL_2$--modular symbols. It would be quite interesting
to extend it to groups of higher rank,
along the lines of [AB] and [AR].

\medskip

{\bf 0.3. Summary of this paper.} I recall the basic
properties of iterated integrals of holomorphic 1--forms
on a simply connected Riemannian surface in \S 1.
The shuffle relations
for the iterated integrals are reflected directly
in terms of a generating function $J$ stating that it is a group--like
element with respect to a comultiplication, cf. Proposition 1.3.1.

\smallskip

Then I turn to the main object of study. 
In \S 2  I define 1--forms of modular and cusp modular type,
introduce and study the iterated and total Mellin transform for
families of such forms. The functional equation for the total Mellin transform is deduced which extends the classical functional equation
for $L$--series.

\smallskip

Using only critical values of these Mellin transforms,
I introduce in 2.6 an iterated modular symbol as a certain noncommutative
1--cohomology class of the relevant subgroup of the
modular group. 

\smallskip

In \S 3, I study the representation of such Mellin transforms
at integer values of their Mellin arguments
in terms of multiple Dirichlet series. The results differ from
the classical ones expressed by the identity (0.1) = (0.2) in two essential
respects. First, iterated integrals are only linear combinations
of certain  multiple Dirichlet series. Second, the latter are {\it not}
of the usual type
$$
\sum_{0<n_1<\dots <n_k}
\frac{a_{1,n_1}\dots a_{n,n_k}}{n_1^{m_1}\dots n_k^{m_k}},
$$
in fact, their coefficients depend on pairwise differences $n_j-n_i$.

\smallskip

In \S  4, the properties of the multiple Dirichlet series
which emerged in \S 3, are axiomatized, and the shuffle relations for them
are deduced.
This requires, however, a considerable extension
of the initial supply of series; the system of those coming
from  1--forms of modular type is not closed.

\smallskip

 \S 5 is dedicated to the iterated analogs of the 
so called Eichler--Shimura and Manin relations for periods of cusp forms.
Whereas the relations of the first type are quite straightforward,
the relations of the second type, involving Hecke operators,
are not obvious. The results presented here (Theorem 5.3) are preliminary, they clearly allow  generalizations and deserve further study.

\smallskip

Finally, in \S 6 I return to the formalism of \S 1 and extend it
by allowing our integrands to have logarithmic singularities
at the boundary. A version of the regularization procedure
I use here is the same as in Drinfeld's paper [Dr2].
It exploits complex analyticity in place of Boutet de Monvel's technique
of [De] and [Ra2]. 

\smallskip

Using the Manin--Drinfeld theorem 
on cusps, I suggest a generalization of Drinfeld's associator
and extend to this case a part of the identities satisfied by the latter.
This list includes the group--like property, the duality, 
and the hexagonal relation, which turn out to have the same
source as the Shimura -- Eichler relations for the periods of cusp forms.
To the contrary, the pentagonal relation seems to be specific 
for the original Drinfeld's associator.  

\medskip

{\it Acknowledgements.} I am grateful to A.~Levin and A.~Goncharov
who have read the first draft of this paper
and made a number of useful suggestions which are
incorporated in the text.

\bigskip

\centerline{\bf \S 1. Iterated integrals of holomorphic 1--forms}

\medskip

{\bf 1.1. Setup.} Let $X$ be a connected Riemann surface,
not necessarily compact, $\Cal{O}_X$ its structure sheaf of holomorphic functions, $\Omega^1_X$ the sheaf of holomorphic 1--forms.
If $\omega$ is a (local) 1--form, $z\in X$ a point,
$\omega (z)$ denotes the value of $\omega$ at $z$, i.e. the
respective cotangent vector.

\smallskip

Let $V$ be a finite set which will be used as a set indexing
various families. Consider
the completed unital semigroup ring freely generated  by $V$. 
We will write it as
the ring of associative formal series $\bold{C}\langle\langle A_V\rangle\rangle$ where $A_V:=(A_v\,|\, v\in V)$ are noncommuting
free formal variables.

\smallskip

More generally, we may consider the ring 
$\Cal{O}_X(U)\langle\langle A_V\rangle\rangle$ where $\Cal{O}_X(U)$
is the ring of holomorphic functions on an open subset $U\subset X$ ($A_v$ commute with 
$\Cal{O}_X(U)$), and the bimodule $\Omega^1_X(U)\langle\langle A_V\rangle\rangle$
over this ring, connected by the differential $d$ such that
$dA_v=0$ for all $v\in V$. Varying $U$, we will get two presheaves;
the sheaves associated with these presheaves are denoted
$\Cal{O}_X\langle\langle A_V\rangle\rangle$, resp. 
$\Omega^1_X\langle\langle A_V\rangle\rangle$, and $d$ extends to them,
so that $\roman{Ker}\,d$ is the constant sheaf 
$\bold{C}\langle\langle A_V\rangle\rangle$.

\smallskip

Let $\omega_V:=(\omega_v\,|\,v\in V)$
be a family of  1--forms holomorphic in $U$ and indexed by $V$. Put
$$
\Omega := \sum_{v\in V} A_v\omega_v\,.
\eqno(1.1)
$$
The total iterated integral of this form
along a piecewise smooth path $\gamma :\,[0,1]\to U$ is
denoted $J_{\gamma}(\Omega )$ or $J_{\gamma}(\omega_V)$ and
 is defined by the following formula:
$$
J_{\gamma}(\Omega):= 1+\sum_{n=1}^{\infty}
\int_{0}^1 \gamma^*(\Omega )(t_1)
\int_{0}^{t_1} \gamma^*(\Omega )(t_2)\dots
\int_{0}^{t_{n-1}}\gamma^*(\Omega )(t_n) \in \bold{C}\langle\langle
A_V\rangle\rangle \,
\eqno(1.2)
$$
where the integration is taken over the simplex $0<t_n<\dots <t_1<1$.
If $\gamma$, $\gamma^{\prime}$  with the same ends are homotopic,
$J_{\gamma}(\Omega)=J_{\gamma^{\prime}}(\Omega)$. 

\smallskip

Putting $z_i=\gamma (t_i)\in X$, $a=\gamma (0)$,
$z=\gamma (1)$, and considering the whole integral
as a function of a variable $z$ we will also write (1.2)
in the form
$$
J_a^z(\Omega )= J_a^z(\omega_V)=1+\sum_{n=1}^{\infty}
\int_{a}^z \Omega (z_1)
\int_{a}^{z_1} \Omega (z_2)\dots
\int_{a}^{z_{n-1}}\Omega (z_n)\,.
\eqno(1.3)
$$
If $U$ is connected and simply connected, this expression is an unambiguously defined element of $\Cal{O}_X(U)\langle\langle A_V\rangle\rangle$.
Otherwise it is a multivalued function of $z$ in this domain.

\smallskip

The following result is classical.

\medskip

{\bf 1.2. Proposition.} {\it (i) $J_a^z(\Omega )$ as a function of $z$
satisfies the equation
$$
dJ_a^z(\Omega )=\Omega (z)\,J_a^z(\Omega ).
\eqno(1.4)
$$
In other words, $J_a^z(\Omega )$ is a horizontal (multi)section of
the flat connection $\nabla_{\Omega} := d -l_{\Omega}$ on
$\Cal{O}_X\langle\langle A_V\rangle\rangle$, where $l_{\Omega}$
is the operator of left multiplication by $\Omega$.

\smallskip

(ii) If $U$ is a simply connected neighborhood of $a$, 
$J_a^z(\Omega )$ is the only horizontal section
with initial 
condition $J_a^a=1$. Any other horizontal section $K^z$ can be 
uniquely written in the form $J_a^z(\Omega )C$, $C\in  
\bold{C}\langle\langle
A_V\rangle\rangle$. In particular, for any $b\in U$,
$$
J_b^z(\Omega )= J_a^z(\Omega )J_b^a(\Omega )
\eqno(1.5)
$$
}
\smallskip

{\bf Proof.} (i) follows directly from (1.3). Since
$J_a^z(\Omega )$ is an invertible element of the ring
$\Cal{O}_X(U)\langle\langle A_V\rangle\rangle$,
we can form $J_a^z(\Omega )^{-1}K^z$ and then directly check
that $d(J_a^z(\Omega )^{-1}K^z )=0$. Hence this element belongs to
$\bold{C}\langle\langle
A_V\rangle\rangle$, and moreover equals its value at $z=a$
that is, $K^a$. Choosing $K^z=J^z_b(\Omega )$, we get (1.5).

\medskip

{\bf 1.3.  $J_a^z(\Omega )$ as a generating series.} Clearly, we have
$$
J_a^z(\omega_V)=J_a^z(\Omega )=1+\sum_{n=1}^{\infty}\sum_{(v_1,\dots ,v_n)\in V^n} 
A_{v_1}\dots A_{v_n}\,
 I_a^z( \,\omega_{v_1},\dots ,\omega_{v_n})\,,
\eqno(1.6)
$$
where
$$
I_a^z( \,\omega_{v_1},\dots ,\omega_{v_n}) =
\int_{a}^{z} \omega_{v_1}(z_1)
\int_{a}^{z_1} \omega_{v_2}(z_2)\dots
\int_{a}^{z_{n-1}}\omega_{v_n}(z_n)
\eqno(1.7)
$$
are the usual iterated integrals.

\smallskip

In the remaining part of this section, and in the main
body of the paper, we will encode various (infinite families of)
relations among the iterated integrals (1.7) in the form
of relations between the generating functions $J_a^z(\omega_V)$.
Generally, our relations between the generating
functions will be (noncommutative) polynomial ones. They may also involve
different families $(\omega_V)$, different integration paths,
and some linear transformations of the formal variables $A_v$:
cf. especially Theorem 2.3, Proposition 5.1.1, Theorem 5.3,
and subsection 6.5 (in the context requiring a regularization).

\medskip

{\bf 1.4. Basic relations between total iterated integrals.}
There are three types of basic relations, which we will call
{\it group--like property, cyclicity}, and {\it functoriality} respectively.

\medskip

{\bf 1.4.1. Proposition.} {\it Consider the comultiplication
$$
\Delta :\, \bold{C}\langle\langle A_V\rangle\rangle\to 
\bold{C}\langle\langle A_V\rangle\rangle\widehat{\otimes}_{\bold{C}}
\bold{C}\langle\langle A_V\rangle\rangle, \ \Delta (A_v)=A_v\otimes 1+
1\otimes A_v 
$$
and extend it to the series
with coefficients $\Cal{O}_X$ and $\Omega^1_X$. Then
$$
\Delta\,(J^z_a(\omega_V)) =J^z_a(\omega_V)\otimes_{\Cal{O}_X} 
J^z_a(\omega_V) \,.
\eqno(1.8)
$$
}
\smallskip

{\bf Proof.} Both sides of (1.8) satisfy the equation $dJ=\Delta (\Omega)J$
and have the initial value 1 at $z=a$.

\smallskip

{\bf NB.} Coefficientwise, (1.8) is a compact version of shuffle relations
for iterated integrals (1.7).

\medskip

{\bf 1.4.2. Cyclicity.} {\it Let $\gamma$ be a closed oriented contractible contour in $U$, $a_1, \dots ,a_n$ points along this contour
(cyclically) ordered compatibly with orientation. Then
$$
J^{a_1}_{a_2}(\Omega )J^{a_2}_{a_3}(\Omega )\dots J^{a_{n-1}}_{a_n}(\Omega )
J^{a_n}_{a_1}(\Omega ) =1.
\eqno(1.9)
$$
}

\smallskip

This follows from (1.5) by induction.

\medskip

{\bf 1.4.3. Functoriality.} Consider an automorphism
$g:\,X\to X$ such that $g^*$ maps into itself the linear space spanned
by $\omega_v$. In particular, there is a constant matrix $G=(g_{vu})$
with rows and columns labeled by $V$ such that
$g^*(\omega_v )=\sum_u g_{vu}\omega_u$. Define
the automorphism $g_*$ of any of the ring/module
of formal series $\bold{C}\langle\langle A_V\rangle\rangle$,
$\bold{C}(X)\langle\langle A_V\rangle\rangle$, $\Omega^1(X)\langle\langle A_V\rangle\rangle$ by the formula $g_*(A_u)=\sum_v A_vg_{vu}$.
On coefficients $g_*$ acts identically.

\medskip

{\bf 1.4.4. Claim.} {\it We have
$$
J_{ga}^{gz}(\omega_V)=g_*(J_a^z(\omega_V))\,.
\eqno(1.10)
$$}

\smallskip

{\bf Proof.} In fact, both sides coincide with $J_a^z(g^*(\omega_V))$.
We will give below a calculation which proves a slightly more general
statement.

\medskip

{\bf 1.5. A variant: multiple lower integration limits.}
Somewhat more generally, in the simply connected case we can
consider a family of points $(a_{\bullet}):=(a_{i,v})$ in $X$ 
indexed by pairs $i=1,2,3, ...;\, v\in V.$

\smallskip

Given such a family and $\omega_V$, we can construct the following
formal series in $\bold{C}(X)\langle\langle A_V\rangle\rangle$
with constant term 1:
$$
J_{(a_{\bullet})}^z(\omega_V):= \sum_{n=0}^{\infty}\sum_{(v_1,\dots ,v_n)\in V^n}
A_{v_1}\dots A_{v_n} I^z_{a_{1,v_1},\dots ,a_{n,v_n}}(\omega_{v_1},\dots ,\omega_{v_n})\,,
\eqno(1.11)
$$
where
$$
I^z_{a_{1,v_1},\dots ,a_{n,v_n}}(\omega_{v_1},\dots ,\omega_{v_n}):=
\int_{a_{1,v_1}}^{z} \omega_{v_1}(z_1)
\int_{a_{2,v_2}}^{z_1} \omega_{v_2}(z_{2})\dots
\int_{a_{n,v_n}}^{z_{n-1}}\omega_{v_n}(z_n)\,.
\eqno(1.12)
$$
As above, $z\in X$ denotes a variable point, the argument of our functions.
Then we have
$$
dJ_{(a_{\bullet})}^z(\omega_V) =\Omega\,J_{(a_{\bullet})}^z(\omega_V)
\eqno(1.13)
$$
and 
$$
J^z_{(a_{\bullet})}(\omega_V)= J^z_a(\omega_V)\,J^a_{(a_{\bullet})}(\omega_V)\,.
\eqno(1.14)
$$
The series (1.11) satisfies the following functoriality
relation generalizing (1.10):

\medskip

{\bf 1.5.1. Claim.} {\it We have
$$
J_{(ga_{\bullet})}^{gz}(\omega_V)=g_*(J_{(a_{\bullet})}^z(\omega_V))\,.
\eqno(1.15)
$$}
 
{\bf Proof.} We will check that both sides coincide
with $J_{(a_{\bullet})}^z(g^*(\omega_V)).$ In fact,
$\int_{gu}^{gv}\nu (z)=\int_{u}^{v}\nu (gz)$ so that,
removing $g$ step by step from the integration limits, we get
$$
I^{gz}_{ga_{1,v_1},\dots ,ga_{n,v_n}}(\omega_{v_1},\dots ,\omega_{v_n})
=I^z_{a_{1,v_1},\dots ,a_{n,v_n}}(g^*(\omega_{v_1}),\dots ,g^*(\omega_{v_n}))\,.
$$
Multiplying the l.h.s. by $A_{v_1}\dots A_{v_n}$ and summing, 
we get the l.h.s. of (1.15).

\smallskip

On the other hand, 
$$
\sum_{v_1,\dots ,v_n\in V^n}A_{v_1}\dots A_{v_n} I^z_{a_{1,v_1},\dots ,a_{n,v_n}}(g^*(\omega_{v_1}),\dots ,g^*(\omega_{v_n}))=
$$
$$
\sum_{v_1,\dots ,v_n\in V^n}A_{v_1}\dots A_{v_n} I^z_{a_{1,v_1},\dots ,a_{n,v_n}}(\sum_{u_1\in V}g_{v_1,u_1}\omega_{u_1},\dots ,
\sum_{u_n\in V}g_{v_n,u_n}\omega_{u_n})=
$$
$$
=\sum \Sb v_1,\dots ,v_n\in V^n\\u_1,\dots ,u_n\in V^n\endSb
A_{v_1}g_{v_1,u_1}\dots A_{v_n}g_{v_n,u_n}
I^z_{a_{1,v_1},\dots ,a_{n,v_n}}(\omega_{v_1},\dots ,\omega_{v_n})
 =
$$
$$
=g_*\left(\sum_{v_1,\dots ,v_n\in V^n}A_{v_1}\dots A_{v_n}
I^z_{a_{1,v_1},\dots ,a_{n,v_n}}(\omega_{v_1},\dots ,\omega_{v_n}) \right)\, .
$$
Summation over $n$ produces the r.h.s. of (1.15), proving the lemma.

\medskip

{\bf 1.6. A variant: nonlinear $\Omega$.}
Let now $\Omega\in \Omega^1_X(U)\langle\langle A_V\rangle\rangle$
be an arbitrary form without a constant term in $A_v$:
$$
\Omega =  \sum_{n=1}^{\infty}\sum_{(v_1,\dots ,v_n)\in V^n}
A_{v_1}\dots A_{v_n}\Omega_{v_1,\dots ,v_n}\,,
\eqno(1.16)
$$
where $\Omega_{v_1,\dots ,v_n}\in \Omega^1_X(U)$.

\smallskip

The total iterated integrals $J_{\gamma}(\Omega )$ and $J_a^z(\Omega )$
are defined by exactly the same formulas (1.2) and (1.3).
It is not true anymore that the coefficients of this series are
the usual iterated integrals. However, an analog of the
Proposition 1.2 and the cyclic identity remain true:

\medskip

{\bf 1.6.1. Proposition.} {\it $J_a^z(\Omega )$ as a function of $z$
satisfies the equation
$$
dJ_a^z(\Omega )=\Omega (z)\,J_a^z(\Omega ).
\eqno(1.17)
$$
If $U$ is a simply connected neighborhood of $a$, 
$J_a^z(\Omega )$ is the only horizontal section
with initial 
condition $J_a^a=1$. Any other horizontal section $K^z$ can be 
uniquely written in the form $J_a^z(\Omega )C$, $C\in  
\bold{C}\langle\langle
A_V\rangle\rangle$. In particular, for any $b\in U$,
$$
J_b^z(\Omega )= J_a^z(\Omega )J_b^a(\Omega )
\eqno(1.18)
$$
}          
\medskip

{\bf 1.6.2. Corollary.} {\it Let $\gamma$ be a closed oriented contractible contour in $U$, $a_1, \dots ,a_n$ points along this contour
(cyclically) ordered compatibly with orientation. Then
$$
J^{a_1}_{a_2}(\Omega )J^{a_2}_{a_3}(\Omega )\dots J^{a_{n-1}}_{a_n}(\Omega )
J^{a_n}_{a_1}(\Omega ) =1.
\eqno(1.19)
$$
}
\medskip

Notice in conclusion that the integral formula 
(0.2) for the multiple zeta values
is not quite covered by the formalism reviewed so far because
the integrands in (0.2) have logarithmic poles at the boundary. We will 
return to this situation in \S 6, to which some readers may prefer to 
turn right away. However, for applications to the integration
of cusp forms in \S 2 -- \S 5  the regular case treated here
suffices.

\bigskip

\centerline{\bf \S 2. 1--forms of modular type, iterated Mellin transform,}

\smallskip

\centerline{\bf and noncommutative modular symbols}

\medskip

{\bf 2.1. Setup.} In this section, $X$ will be the upper
half plane $H$ and $z$  the standard complex
coordinate. $H$ is endowed with the metric of constant curvature
$-1$: $ds^2=|dz|^2/(\roman{Im}\,z)^2$.

\smallskip

The limits of integration in our iterated integrals generally lie in $H$,
but may be
``improper'' as well, that is, belong to the set of cusps $\bold{Q}\cup \{i\infty\}$.
If this is the case, we always assume that the respective integration 
path in some neighborhood of the cusp coincides with  a segment of a geodesic curve.

\smallskip

Our 1--forms  generally will have the following 
structure.

\smallskip

{\bf 2.1.1. Definition.} {\it (i) A 1--form $\omega$ on $H$
is called a form of modular type, if it can be represented as
$f(z) z^{s-1}dz$
where $s$ is a complex number, and $f(z)$ is a modular form of some weight
with respect to a congruence subgroup of the modular group.

The modular form $f(z)$ is then well defined and called
the associated modular form (to $\omega$), and the number $s$
is called the Mellin argument of $\omega$.

\smallskip

(ii)   $\omega$ is called a form of cusp modular type 
if the associated $f(z)$ is a cusp form.}

\smallskip

To fix notation, we will recall below some classical facts.

\medskip

{\bf 2.1.2. Action of automorphisms.} Any matrix 
$\gamma \in GL_2^+(\bold{R})$ defines
a holomorphic isometry of $H$, namely $z\mapsto [\gamma\, ]z$ where
$[\gamma ]$ is the fractional linear transformation 
corresponding to $\gamma$. We will denote this automorphism
also $\gamma$. It induces the inverse image maps on
the sheaves $(\Omega^1_H)^{\otimes r}$ of holomorphic tensor differentials
of degree $r$:
$$
\gamma^*(f(z)\, (dz)^{r}) = f([\gamma\, ]z)\,(d [\gamma\, ]z)^r=
(\roman{det}\,\gamma )^r  f([\gamma\, ]z)\,\frac{({dz})^r}{(c_{\gamma}z+
d_{\gamma})^{2r}}
\eqno(2.1) 
$$
where $(c_{\gamma},d_{\gamma})$ is the lower row of $\gamma$. 

\smallskip

If one identifies $(\Omega^1_H)^{\otimes r}$ with $\Cal{O}_H$
by sending $(dz)^{r}$ to 1,  (2.1) turns into the action
of weight $2r$ on functions which is traditionally written as a right action:
$$
f|[\gamma\,]_{2r}(z) :=  
(\roman{det}\,\gamma )^r  f([\gamma\, ]z)\,(c_{\gamma}z+d_{\gamma})^{-2r} .
\eqno(2.2)
$$

Assume that $f(z)(dz)^r$ is invariant with respect to
$\gamma$. Then, writing $f(z)z^{s-1}dz=f(z) (dz)^r\cdot z^{s-1}(dz)^{1-r}$,
we see that 
$$
\gamma^*(f(z)\,z^{s-1} dz)= 
(\roman{det}\,\gamma )^{1-r}f(z) 
(a_{\gamma}z+b_{\gamma})^{s-1}(c_{\gamma}z+d_{\gamma})^{2r-1-s} dz .
\eqno(2.3)
$$
where $(a_{\gamma}, b_{\gamma})$ is the upper row of $\gamma$.
In particular, if $2r\ge 2$ is an integer,
$\gamma^*$ maps into itself the space of 1--forms spanned by
$$
f(z)\,z^{s-1} dz, \quad 1\le s\le 2r-1,\ s\in \bold{Z} .
\eqno(2.4)
$$
More generally, if 
$$
\gamma^*(f(z)(dz)^r)=\chi (\gamma )\,f(z)(dz)^r
\eqno(2.5)
$$
for some $\chi (\gamma )\in \bold{C}$, then 
$$
\gamma^*(f(z)\,z^{s-1} dz)= 
(\roman{det}\,\gamma )^{1-r}f(z) 
\chi (\gamma )\,(a_{\gamma}z+b_{\gamma})^{s-1}(c_{\gamma}z+d_{\gamma})^{2r-1-s} dz ,
\eqno(2.6)
$$
and the space (2.4) will still remain invariant.

\smallskip

We can apply this  formalism to the spaces of modular forms
of weight $2r$
with respect to a congruence subgroup $\Gamma$ of $SL_2(\bold{Z}),$
i.e. to the functions $f$ in $(dz)^{-r}((\Omega^1_H)^{\otimes r})^{\Gamma}.$
Two special cases will be of particular interest:

\medskip

(i) For any such $f$, the space of 1--forms spanned by (2.4) 
is $\Gamma$--invariant.

\smallskip

(ii) Assume that $\Gamma =\Gamma_0(N)$. This group is normalized by the involution
$$
g=g_N:=\left(\matrix 0& -1\\ N& 0\endmatrix \right)
\eqno(2.7)
$$
Therefore this involution maps into itself the space of
$\Gamma_0(N)$--modular forms, and the latter has a basis consisting
of forms with
$$
g_N^*(f(z)(dz)^r)=\varepsilon_f\,f(z)(dz)^r,\quad  \varepsilon_f=\pm 1.
\eqno(2.8)
$$
Applying (2.6) with $\gamma =g_N$ we get, for any complex $s$,
$$
g_N^*(f(z) z^{s-1}dz)=\varepsilon_f N^{r-s} f(z) z^{2r-1-s} dz.
\eqno(2.9)
$$
\medskip

{\bf 2.1.3. Geodesics and cusp forms.} The geodesic
from $0$ to $i\infty$ is the upper half of the pure imaginary line.
The unoriented distance of a point $iy$ on it to $i$
is $|\roman{log}\,y|.$ The exponential of this distance is thus
$y$, if $y>1$, and $y^{-1}$, if $y<1$. If we replace $i$ by another reference point, even outside  of imaginary axis, the exponential of the distance will behave like
$e^{O(1)}y$ (resp. $e^{O(1)}y^{-1}$) as $y\to \infty$ (resp. $y\to 0$.)

\smallskip

Let $f(z)$ be a cusp form of weight $2r$ for a congruence subgroup.
Then it can be represented by a Fourier series
$f(z)=\sum_{n=1}^{\infty}c_ne^{2\pi inz/N}$ for some $N\in \bold{Z}_+$,
whose coefficients are polynomially bounded: $c_n=O(n^C)$ for some 
$C>0.$ Therefore we have $|f(iy)|=O(e^{-ay})$ for some $a>0$
as $y\to \infty$. From the previous analysis
it follows that more generally, for any cusp form
and any geodesic connecting two cusps, $|f(z)|=O(e^{-ay(z)})$
for some $a>0$ as $z$ tends  along the geodesics to one of its ends,
where this time $y(z)$ means the exponentiated geodesic distance from
$z$ to any reference point in $H$, fixed once and for all.

\smallskip

Let now $\omega(z) =f(z)z^{s-1}dz$ be an 1--form of cusp modular type. 
Then the estimates above show
that the following expected properties indeed hold.

\smallskip

a) As $z_0\to i\infty$ along the imaginary axis, the family $\int_{z_0}^{z} \omega$
of holomorphic functions of
$z$ in any  bounded domain $H$ converges absolutely and uniformly
to a holomorphic function of $z$ which is denoted
$\int_{i\infty}^z \omega$. The same remains true, if one replaces $i\infty$
by $0$.

\smallskip

These integrals are holomorphic functions of the Mellin argument
$s$ of $\omega$ as well.

\smallskip

b) The sum  $(\int_{i\infty}^z + \int^{0}_z)\, \omega$ does not depend
on $z$ in $H$ and is denoted $\int_{i\infty}^0 \omega$.
As a function of $s$, it is called {\it the 
classical Mellin transform of $\omega$.}
 
\smallskip

Denote this classical transform $\Lambda (f;s).$ Assume that
$f$ satisfies (2.8). Then we have the classical functional equation
$$
\Lambda (f;s) = -\varepsilon_f N^{r-s} \Lambda (f; 2r-s)\,,
\eqno(2.10)
$$
because in view of (2.9)
$$
\int_{i\infty}^0 \omega =-\int^{i\infty}_0 \omega =
-\int_{g_N(i\infty)}^{g_N(0)} \omega =
-\int_{i\infty}^0 g_N^*(\omega )= -\varepsilon_f N^{r-s}
\int_{i\infty}^0 f(z) z^{2r-1-s} dz.
$$
Another identity in the same vein uses the fact that
$i/\sqrt{N}$ is the fixed point of $g_N$ so that
$\Lambda (f;s)$ can be written as
$$
\Lambda (f;s) = \int_{i\infty}^{i/\sqrt{N}} \omega -
\int_{i\infty}^{i/\sqrt{N}} g_N^*(\omega )\, .
\eqno(2.11)
$$
This allows one to use the Fourier expansions of
$f(z)$ and $f|[g_N]_{2r}(z)$ in order to deduce
series expansions for $\Lambda (f;s)$ (notice that
the Fourier expansions cannot be term--wise integrated near $z=0$ because 
the formal integration produces a divergent series).

\smallskip

Now we can finally write down the analogs of these definitions and 
results for iterated integrals.

\medskip

{\bf 2.2. Definition.} {\it (i) Let $f_1,\dots ,f_k$ be a finite sequence
of cusp forms with respect to a congruence subgroup,
$\omega_j(z):= f_j(z)\, z^{s_j-1} dz.$ The iterated Mellin transform 
of $(f_j)$ is, by definition,
$$
M (f_1,\dots ,f_k;s_1,\dots ,s_k) :=
I_{i\infty}^0(\omega_1,\dots ,\omega_k)=
$$
$$
=\int_{i\infty}^{0} \omega_{1}(z_1)
\int_{i\infty}^{z_{1}} \omega_{2}(z_{2})\dots
\int_{i\infty}^{z_{n-1}}\omega_n(z_n)
\eqno(2.12)
$$
\smallskip

(ii) Let $f_V=(f_v\,|\,v\in V)$ be a finite family
of cusp forms with respect to a congruence subgroup,
$s_V=(s_v\,|\,v\in V)$ a finite family of complex numbers, $\omega_V=(\omega_v)$, where
$\omega_v(z):= f_v(z)\, z^{s_v-1} dz.$
The total Mellin transform 
of $f_V$ is, by definition,
$$
TM(f_V;s_V) := J_{i\infty}^0(\omega_V) =
$$
$$
=\sum_{n=0}^{\infty}\sum_{(v_1,\dots ,v_n)\in V^n}
A_{v_1}\dots A_{v_n}\,M (f_{v_1},\dots ,f_{v_n};s_{v_1},\dots ,s_{v_n})
\eqno(2.13)
$$
(cf. (1.3)).}

\smallskip

Below we will assume that the space  spanned by all $\omega_v$
is stable with respect to some $g_N^*$. Then as in 1.4.3
denote by $G=(g_{vu})$ the matrix of this action on
$(\omega_v)$, and by $g_{N*}$ the action of the transposed matrix
on the formal variables $(A_v)$. 

\smallskip

For example, if $(\omega_v)$ and $(f_v(dz)^{r_v})$ respectively
can be represented as a union of pairs of forms,
corresponding to the left and right hand sides of (2.9),
the matrix $G$ consists of two by two antidiagonal blocks each of which
after the classical Mellin transform produces
a functional equation of the form (2.10).

\medskip

{\bf 2.3. Theorem.} {\it (i) If the space  spanned by all $\omega_v$
is stable with respect to some $g_N^*$, we have
the following functional equation: 
$$
J_{i\infty}^0(\omega_V) =g_{N*}(J_{i\infty}^0(\omega_V))^{-1}\,.
\eqno(2.14)
$$

(ii) In the assumptions of the Definition 2.2 (ii),
denote the weight of $f_v$ by $2r_v$ and assume that $f_v$
is an eigenvector
for $g_N^*$ with eigenvalue $\varepsilon_v$. Then
the total Mellin transform (2.13) satisfies
$$
TM(f_V;s_V)= g_*(TM(f_V;2r_V-s_V))^{-1}
\eqno(2.15)
$$
where $g_*$ multiplies each $A_v$ by $\varepsilon_v N^{r_v-s_v}$.
}
\smallskip

{\bf Proof.} This is a straightforward corollary of the definitions
and formulas (1.9) and (1.10) as soon as one has checked
that the latter formulas are applicable to the improper iterated integrals
of the 1--forms of cusp modular type. 

\smallskip

This check is a routine matter, since at each step
of an iterated integration we multiply the result of the
previous step by a holomorphic function of the type
$f(z) \,z^{s-1}$ which is bounded by $O(e^{-ay(z)})$
as in 2.1.3 above as $z$ tends to $0$ or $i\infty$.

\smallskip

Notice in conclusion that no analog of the functional equation
(2.11) can be written for the individual Mellin transforms (2.12),
because applying $g_N$ to the integration limits in them we get an expression
which is not a Mellin transform in our sense. Only putting them all together
produces the necessary environment for replacing the overall minus sign 
at the r.h.s. of (2.10) by the overall exponent $-1$
at the r.h.s. of (2.15).

\smallskip

A similar reasoning establishes the iterated analog of (2.11):

\medskip

{\bf 2.4. Proposition.} {\it We have
$$
TM(f_V;s_V)= 
(g_{N*}J_{i\infty}^{i/\sqrt{N}} (\omega_V))^{-1}
J_{i\infty}^{i/\sqrt{N}} (\omega_V)\, .
\eqno(2.16)
$$
} 
\medskip

{\bf 2.5. Pushing down iterated integrals.} Let $\omega$ be
an 1--form of modular type whose associated modular
form has weight 2 with respect to a subgroup $\Gamma$ of the modular
group, and whose Mellin argument is 1. In this case
$\omega$ is $\Gamma$--invariant so that it can be pushed down to
an 1--form $\nu$ on $X_{\Gamma}^{\circ}:=\Gamma \setminus H$. Instead of integrating
$\omega$ along a path in $H$, we can integrate $\nu$
along the push--down of this path to $X_{\Gamma}^{\circ}$. If
all $\omega_v$ have this property, all relevant iterated integrals can
be pushed down to $X_{\Gamma}^{\circ}$. 

\smallskip

This argument admits a partial generalization to higher weights.
Assume that the modular form associated with $\omega$
has weight $2r>2$, whereas its Mellin argument
is an integer belonging to the critical strip (2.4),
$1\le s\le 2r-1$. In this case the relevant
simple integral along, say, $\{i\infty ,0\}$
can be pushed down to the Kuga--Sato variety $X_{\Gamma}^{(2r-2)}$
which is the $(2r-2)$--th fibered power of the universal elliptic
curve over $X_{\Gamma}$ or rather its compactified smooth model.
However, on  $X_{\Gamma}^{(2r-2)}$ we obtain an integral
of a holomorphic form $\widehat{\omega}$ of degree $2r-1$ over a relative cycle of the same dimension which is $>1$. Therefore iterated 
``line'' integrals of such forms
on $H$ cannot be directly translated into integrals of the same type
on $X_{\Gamma}^{(2r-2)}$. 

\smallskip

On the other hand, one can generally define Chen's iterated integrals of forms 
of arbitrary degree, say, $\widehat{\omega}_v$ on $X_{\Gamma}^{(2r-2)}$,
which take values in the space of differential forms on the
path space $PX_{\Gamma}^{(2r-2)}$ and not just $\bold{C}$: cf.
reports [Ch] and [Ha], as well as references therein.
Studying properties of such iterated integrals in the modular case like presents
an interesting challenge.

\smallskip

Here I will restrict myself to explaining how $\widehat{\omega}$
looks like and why its periods coincide with integrals of $\omega$
along geodesics. 
For more details, see [Sh1], [Sh2], and especially [Sh3].

\smallskip

Denote by $\Gamma^{(r)}$ the semidirect product
$\Gamma\ltimes (\bold{Z}^{2r-2}\times \bold{Z}^{2r-2})$
acting upon $H\times \bold{C}^{2r-2}$ via
$$
(\gamma ;\,n,m)\,(z,\zeta ):= ([\gamma ]z;\,(c_{\gamma}z +d_{\gamma})^{-1}
(\zeta +zn +m))\,.
$$
Here $n=(n_1,\dots ,n_{2r-2})$, 
$m=(m_1,\dots ,m_{2r-2})$, $\zeta =(\zeta_1,\dots ,\zeta_{2r-2})$,
and $nz=(n_1z,\dots ,n_{2r-2}z)$.

\smallskip

If $f(z)$ is a holomorphic modular form of weight $2r$, then
$f(z)dz\wedge d\zeta_1\wedge \dots \wedge d\zeta_{2r-2}$
is a $\Gamma^{(r)}$--invariant holomorphic volume form
on $H\times \bold{C}^{2r-2}$. Hence one can push it down
to (a Zariski open subset of) the quotient
$\Gamma^{(r)} \setminus (H\times \bold{C}^{2r-2})$ which
is a Zariski open subset of the respective Kuga--Sato variety.
Denote by $\widehat{\omega}$ the image of this form.
Notice that it is common for all 1--forms of modular type
$\omega =f(z)z^{s-1}dz$ with different Mellin arguments $s$.

\smallskip

A detailed analysis of singularities performed in [Sh2], [Sh3]
shows that the map $f\mapsto \widehat{\omega}$ induces an isomorphism
of the space of cusp forms  of weight $2r$ with the space
of holomorphic volume forms on an appropriate smooth
projective Kuga--Sato variety. (As I have already remarked in the 
Introduction, it would be useful to replace it
by the base extension  $(\overline{M}_{1, 2r-2})_{X_{\Gamma}}$.)

\smallskip

The dependence of the period of $\omega$ on the integration
path and on the Mellin argument is reflected in the choice of
the relative cycle over which we integrate $\widehat{\omega}$. 

\smallskip

More precisely, let $\alpha ,\beta \in \bold{P}^1(\bold{Q})$ 
be two cusps in $\overline{H}$ and let $p$ be a geodesic joining
$\alpha$ to $\beta$. Fix $(n_i)$ and $(m_i)$ as above.
Construct a cubic singular cell $p\times (0,1)^{2r-2}\to H\times\bold{C}^{2r-2}$:
$(z, (t_i))\mapsto (z, (t_i(zn_i+m_i)))$. Take the
$S_{2r-2}$--symmetrization of this cell and
push down the result to the Kuga--Sato variety.
We will get a relative cycle whose homology class is Shokurov's
higher modular symbol $\{\,\alpha, \beta ;\,n,m\}_{\Gamma}.$
From this construction it is almost obvious that
$$
\int_{\alpha}^{\beta} f(z) \sum_{i=1}^{2r-2}(n_iz+m_i)\,dz=
\int_{\{\,\alpha, \beta ;\,n,m\}_{\Gamma}} \widehat{\omega}\,.
$$
The singular cube $(0,1)^{2r-2}$ may also be replaced by
an evident singular simplex. This can be useful for transposing the results
of [GoMa] to the genus one moduli spaces. 

\medskip

{\bf 2.6. Noncommutative modular symbols and continued fractions.} 
I will define in this subsection a generalization of
modular symbols involving iterated integrals and allowing
a mixture of forms of different weights with respect to the same
subgroup $\Gamma$ of $SL(2,\bold{Z})$. 

\smallskip

Let $(\omega_v)$ be a family of linearly independent
1--forms of cusp modular type
whose Mellin arguments are integers lying in the respective
critical strip as in (2.4). Let $\Gamma$ be a subgroup of
modular group acting on the space spanned by $(\omega_v)$
as in (2.3). Denote by $\Pi$ the multiplicative group
of power series in $(A_v)$ with constant term 1. Clearly, the map
$J\mapsto g_*J$ (see 1.4.3.) defines the left action of $\Gamma$ 
on $\Pi$. 

\medskip

{\bf 2.6.1. Proposition--Definition.} {\it
(i) For each $a\in \bold{P}^1(\bold{Q})$, the map
$\Gamma\to \Pi :\, \gamma\mapsto J^a_{\gamma a}(\Omega )$ is a noncommutative
1--cocycle $\zeta_a$ in $Z^1(\Gamma ,\Pi )$.

\smallskip

(ii) The cohomology class of $\zeta_a$ in $H^1(\Gamma ,\Pi )$
does not depend on the choice of 
$a$ and is called the noncommutative modular symbol.}

\smallskip

{\bf Proof.} We have, omitting $\Omega$ for brevity,
and using (1.9), (1.10):
$$
J^a_{\gamma\beta a} =J^a_{\gamma a} J^{\gamma a}_{\gamma \beta a} = J^a_{\gamma a} \gamma_*(J^a_{\beta a})
$$
which means that $\zeta_a$ is an 1--cocycle. Moreover, if $b$ is
another cusp,
$$
J^a_{\gamma a} = J^a_{b}J^b_{\gamma b} J^{\gamma b}_{\gamma a} = J^a_b J^b_{\gamma b} (\gamma_*(J^a_{b}))^{-1}
$$
that is, $\zeta_a$ and $\zeta_b$ are homologous.

\smallskip

{\bf Remark.} Assume that the cusp forms associated with
$(\omega_v)$ span the sum of all spaces of cusp forms of certain weights,
and for each weight and each cusp form, all admissible Mellin arguments actually occur. Then the linear in $A_v$ term of $\zeta_a$
encodes all periods of the involved cusp forms
along all classical modular symbols corresponding to loops in $X_{\Gamma}(\bold{C})$
starting and ending at the cusp $\Gamma a$.

\medskip

{\bf 2.6.2. Iterated integrals between arbitrary cusps.}
The group $\Gamma$  generally does not act transitively on cusps,
so that the components of cocycles $\zeta_a$ do not contain
iterated integrals along all geodesics connecting two cusps.
One can use the technique of continued fractions
as in [Ma1], [Ma2] in order to express all such integrals 
through a finite number of them.

\smallskip

Namely, choose a set of representatives $C$ of left cosets
$\Gamma \setminus SL_2(\bold{Z})$. Call the iterated integrals
of the form $(J_{g(i\infty )}^{g(0)})^{\pm 1}$, $g\in C$, {\it primitive} ones.
Notice that when $g\notin \Gamma$ the space spanned by
$(\omega_v)$ is not generally $g^*$--stable so that we cannot define
$g_*$.

\medskip

{\bf 2.6.3. Proposition.} {\it Each $J^a_b$ can be expressed as
a noncommutative monomial in $\gamma_*(J^c_d)$ where
$\gamma$ runs over $\Gamma$ and $J^c_d$ runs over primitive integrals.}

\smallskip

{\bf Proof.} First, we can write $J^a_b=(J^{i\infty}_a)^{-1}J^{i\infty}_b$.
So it remains to find a required expression for $J_{i\infty}^a$.
Assume that $a>0$; the case $a<0$ can be treated similarly.
Consider the consequent convergents to $a$:
$$
a=\frac{p_n}{q_n},\ \frac{p_{n-1}}{q_{n-1}},\ \dots ,\ \frac{p_0}{q_0}=
\frac{p_0}{1},\ \frac{p_{-1}}{q_{-1}}:=\frac{1}{0}.
$$
Put
$$
g_k:=\left(\matrix p_{k} & (-1)^{k-1}p_{k-1}\\q_{k} & (-1)^{k-1}q_{k-1} \endmatrix
\right),\quad k=0, \dots , n.
$$
We have $g_k=g_k(a) \in SL_2(2,\bold{Z}).$
Put $g_k=\gamma_kc_k$ where $\gamma_k\in\Gamma$ and
$c_k\in C$ are two sequences of matrices depending on $a$.
Then from cyclicity we get
$$
J_{i\infty}^a=\prod_{k=0}^nJ_{p_{k-1}/q_{k-1}}^{p_k/q_k}
$$
and from functoriality we obtain
$$
J_{p_{k-1}/q_{k-1}}^{p_k/q_k} = \gamma_{k*}(J^{c_k(i\infty)}_{c_k(0)} ).
$$

\bigskip

\centerline{\bf \S 3. Values of iterated Mellin transforms
at integer points}

\smallskip

\centerline{\bf and multiple Dirichlet series}

\medskip

In this section, we collect some formulas expressing 
iterated Mellin transforms (2.12) at integer values
of their Mellin arguments as  linear combinations
of ``multiple Dirichlet series''.  

\medskip

{\bf 3.1. Notation.}
Consider a family of 1--forms $\omega_V$, $v\in V$,
satisfying the following conditions. First,
$$
\omega_v(z)=\sum_{n=1}^{\infty} c_{v,n}e^{2\pi inz} z^{m_v-1}dz,\quad 
c_{v,n}\in \bold{C}\,,\quad m_v\in \bold{Z},\, m_v\ge 1.
\eqno(3.1)
$$
Moreover, assume that $c_{v,n}=O(n^{C})$ for some $C$ and each $v$.

\smallskip

Until a problem of analytic continuation arises, we do not have to assume
modularity. The notation $m_v$ replacing former $s_v$ 
is chosen to remind that these Mellin arguments are natural
numbers.

\smallskip

We start with introducing some notation.

\medskip

{\bf 3.1.1. Functions $L(z;\omega_{v_k},\dots ,\omega_{v_1};
j_k,\dots ,j_1)$.} Choose $k\ge 1;\,v_k,\dots ,v_1\in V$,
and nonnegative integers $j_k,\dots ,j_1$; it is convenient 
to add $j_0=0.$ In our applications,
$j_a$ will satisfy the following restrictions:
$$
j_a\le m_{v_a}-1+j_{a-1}.
\eqno(3.2)
$$

\smallskip

Now put
$$
L(z;\omega_{v_k},\dots ,\omega_{v_1};
j_k,\dots ,j_1):=
$$
$$
= (2\pi iz)^{j_k} \sum_{n_1,\dots ,n_k\ge 1}
\frac{c_{v_1,n_1}\dots c_{v_k,n_k}e^{2\pi i(n_1+\dots +n_k)z}}{n_1^{m_{v_1}+j_0-j_1}(n_1+n_2)^{m_{v_2}+j_1-j_2} \dots (n_1+\dots +n_k)^{m_{v_k}+j_{k-1}-j_k}}
\,.
\eqno(3.3)
$$
Thanks to the presence of exponential terms in (3.2), this series
absolutely converges for any $z$ with $\roman{Im}\,z >0$
and defines a holomorphic function in $H$.

\smallskip

Notice that the enumeration of arguments of $L$ is reversed
in order to get a more natural enumeration of factors
in the summands of (3.3).

\medskip

{\bf 3.1.2. Numbers $L(0;\omega_{v_k},\dots ,\omega_{v_1}; j_k,
j_{k-1},\dots ,j_1)$.} If we formally put $z=0$ in 
the expansion for 
$$
 (2\pi iz)^{-j_k}\,L(z;\omega_{v_k},\dots ,\omega_{v_1};
j_k,\dots ,j_1)
$$
we will get the formal series
$$
 \sum_{n_1,\dots ,n_k\ge 1}
\frac{c_{v_1,n_1}\dots c_{v_k,n_k}} {n_1^{m_{v_1}+j_0-j_1}(n_1+n_2)^{m_{v_2}+j_1-j_2} \dots 
(n_1+\dots +n_k)^{m_{v_k}+j_{k-1}-j_k}}\,.
\eqno(3.4)
$$
We have
$$
c_{v_1,n_1}\dots c_{v_k,n_k}= O((n_1n_2\dots n_k)^C)=
O((n_1+\dots +n_k)^{kC})\,.
\eqno(3.5)
$$
Assume that (3.2) holds. Then  the
general term of (3.4) is bounded by
$$
\frac{1}{n_1(n_1+n_2)\dots (n_1+\dots +n_{k-1})
(n_1+\dots +n_k)^{m_{v_k}+j_{k-1}-j_k-1-kC}}.
$$
Hence (3.4) absolutely converges as long as
$$
m_{v_k}+j_{k-1} -j_k > 1+kC\,.
$$
Summarizing, we get three alternatives, describing the 
possible behavior of 
$$
L(z;\omega_{v_k},\dots ,\omega_{v_1}; j_k,
j_{k-1},\dots ,j_1)
$$ 
as $z\to 0$. We will later
identify the respective limit as a component of the total
Mellin transform of $(\omega_V)$:

\smallskip

{\it Case 1: $j_k=0$ and $m_{v_k}+j_{k-1} > 1+kC$.}
Then the limit exists, and equals to the ``multiple Dirichlet series'' 
$$
L(0;\omega_{v_k},\dots ,\omega_{v_1}; 0,
j_{k-1},\dots ,j_1) =
$$
$$
=  \sum_{n_1,\dots ,n_k\ge 1}
\frac{c_{v_1,n_1}\dots c_{v_k,n_k}} {n_1^{m_{v_1}+j_0-j_1}(n_1+n_2)^{m_{v_2}+j_1-j_2} \dots 
(n_1+\dots +n_k)^{m_{v_k}+j_{k-1}-j_k}}\,.
\eqno(3.6)
$$
($j_0=j_k=0$ appear at the r.h.s. only for the uniformity of
notation).

\smallskip

{\it Case 2: $j_k>0$ and $m_{v_k}+j_{k-1} -j_k > 1+kC$.}
Then the limit exists, and vanishes thanks to the factor
$(2\pi iz)^{j_k}$: 
$$
L(0;\omega_{v_k},\dots ,\omega_{v_1}; j_k,
j_{k-1},\dots ,j_1) =0
\eqno(3.7)
$$

\smallskip

{\it Case 3: $j_k>0$ and $m_{v_k}+j_{k-1} -j_k \le 1+kC$.}
In this case an additional study is needed. 

\medskip

\smallskip

We can now formulate the first main result of this section.

\medskip

{\bf 3.2. Theorem.} {\it For any $k\ge 1$, $(v_1,\dots ,v_k)\in V^k$,
and $\roman{Im}\,z>0$ we have
$$
(2\pi i)^{m_{v_1}+\dots +m_{v_k}}I_{i\infty}^z(\omega_{v_k},\dots ,\omega_{v_1})=
$$
$$
=
(-1)^{\sum_{i=1}^k (m_{v_i}-1)}
\sum_{j_1=0}^{m_{v_1}-1} \sum_{j_2=0}^{m_{v_2}-1+j_1}\dots
\sum_{j_k=0}^{m_{v_k}-1+j_{k-1}}
(-1)^{j_k}\, \times
$$
$$
\times \, \frac{(m_{v_1}-1)!(m_{v_2}-1+j_1)!\dots (m_{v_k}-1+j_{k-1})!}{j_1!
j_2! \dots j_k!}\,L(z;\omega_{v_k},\dots ,\omega_{v_1};
j_k,\dots ,j_1)\,.
\eqno(3.8)
$$}

\smallskip

The proof requires an auxiliary construction.

\medskip

{\bf 3.3. Auxiliary polynomials $D_{m_{v_1},\dots ,m_{v_k}}^{n_1,\dots ,n_k}(t)$.} Choose now $k\ge 1;\,v_k,\dots ,v_1\in V$,
and positive integers $n_k,\dots ,n_1$. It is convenient to
agree that for $k=0$ the respective families are empty.

\smallskip

Define inductively polynomials 
$$
D_{m_{v_1},\dots ,m_{v_k}}^{n_1,\dots ,n_k}(t)\in \bold{Q}[t]
$$
putting 
$D_{\emptyset}^{\emptyset}=1$, and
$$
D_{m_{v_1},\dots ,m_{v_{k+1}}}^{n_1,\dots ,n_{k+1}}(t)= (1+\partial_t)^{-1} 
(D_{m_{v_1},\dots ,m_{v_k}}^{n_1,\dots ,n_k}
\left(\frac{n_1+\dots +n_k}{n_1+\dots +n_{k+1}}\, t\right)\cdot t^{m_{v_{k+1}}-1})\, ,
\eqno(3.9)
$$
where 
$$
(1+\partial_t)^{-1} :=\sum_{k\ge 0} (-1)^k\partial_t^k
$$
as a linear operator on polynomials.

\smallskip

For example,
$$
D_{m_{v_1}}^{n_1}(t)=(-1)^{m_{v_1}-1}(m_{v_1}-1)!
\sum_{j_1=0}^{m_{v_1}-1}\frac{(-1)^{j_1}t^{j_1}}{j_1!}\,.
\eqno(3.10)
$$
In particular, $D_{m_{v_1}}^{n_1}(0)=(-1)^{m_{v_1}-1}(m_{v_1}-1)!$
Furthermore, 
$$
D_{m_{v_1},m_{v_2}}^{n_1,n_2}(t)=
(-1)^{m_{v_1}-1}(m_{v_1}-1)!
\sum_{j_1=0}^{m_{v_1}-1}\frac{(-1)^{j_1}}{j_1!}
(1+\partial_t)^{-1}\,\frac{n_1^{j_1}t^{j_1+m_{v_2}-1}}{(n_1+n_2)^{j_1} }
=
$$
$$
=(-1)^{m_{v_1}-1+m_{v_2}-1}(m_{v_1}-1)! (m_{v_2}-1)!
\sum_{j_1=0}^{m_{v_1}-1}\frac{(-1)^{j_1}n_1^{j_1}}{j_1!
(n_1+n_2)^{j_1}}\sum_{j_2=0}^{m_{v_2}-1+j_1}
\frac{(-1)^{j_2}t^{j_2}}{j_2!}
\eqno(3.11)
$$
In particular,
$$
D_{m_{v_1},m_{v_2}}^{n_1,n_2}(0) =
(-1)^{m_{v_1}-1+m_{v_2}-1}(m_{v_1}-1)! (m_{v_2}-1)!
\sum_{j_1=0}^{m_{v_1}-1}\frac{(-1)^{j_1}n_1^{j_1}}{j_1!
(n_1+n_2)^{j_1}}
$$

The general formula looks as follows:

\medskip

{\bf 3.3.1. Proposition.} {\it We have for $k\ge 1$:
$$
D_{m_{v_1},\dots ,m_{v_{k}}}^{n_1,\dots ,n_{k}}(t)=
(-1)^{\sum_{i=1}^k (m_{v_i}-1)}
\sum_{j_1=0}^{m_{v_1}-1} \sum_{j_2=0}^{m_{v_2}-1+j_1}\dots
\sum_{j_k=0}^{m_{v_k}-1+j_{k-1}}
(-1)^{j_k}\, \times
$$
$$
\times \, \frac{(m_{v_1}-1)!(m_{v_2}-1+j_1)!\dots (m_{v_k}-1+j_{k-1})!}{j_1!
j_2! \dots j_k!}\,\times
$$
$$
\times \,
\frac{1}{n_1^{-j_1}(n_1+n_2)^{j_1-j_2}\dots (n_1+\dots +n_{k-1})^{j_{k-2}-j_{k-1}}
(n_1+\dots +n_k)^{j_{k-1}}}\,t^{j_k}\,.
\eqno (3.12)
$$
}
{\bf Proof.} We argue by induction on $k$. Assume that (3.12) 
holds for $k$ and apply
the operator at the right hand side of (3.9) to the right hand side
of (3.12). Looking for brevity only at the last line of (3.11),
we get:
$$
\frac{1}{n_1^{-j_1}(n_1+n_2)^{j_1-j_2}\dots 
(n_1+\dots +n_k)^{j_{k-1}-j_k}
(n_1+\dots +n_{k+1})^{j_k}}\,(1+\partial_t)^{-1}t^{m_{v_{k+1}}-1+j_k}=
$$
$$
=\frac{1}{n_1^{-j_1}(n_1+n_2)^{j_1-j_2}\dots 
(n_1+\dots +n_k)^{j_{k-1}-j_k}
(n_1+\dots +n_{k+1})^{j_k}}\,\times
$$
$$
\times \,\sum_{j_{k+1}=0}^{m_{v_{k+1}}-1+j_k}
\frac{(-1)^{j_{k+1}}(-1)^{m_{v_{k+1}}-1+j_k}(m_{v_{k+1}}-1+j_k)!}{j_{k+1}!} t^{j_{k+1}}\,.
$$
Combining this with (3.12) for $k$, we get (3.12) for $k+1$.

\smallskip

{\bf 3.4. Proof of Theorem 3.2.} By induction on $k$ we will prove the following
formula:
$$
(2\pi i)^{m_{v_1}+\dots +m_{v_k}}I^z_{i\infty}(\omega_{v_k},\dots ,\omega_{v_1})=
$$
$$
=\sum_{n_1,\dots ,n_k\ge 1}
\frac{c_{v_1,n_1}\dots c_{v_k,n_k}e^{2\pi i(n_1+\dots +n_k)z}}{n_1^{m_{v_2}}(n_1+n_2)^{m_{v_1}} \dots (n_1+\dots +n_k)^{m_{v_k}}}
\, D_{m_{v_1},\dots ,m_{v_k}}^{n_1,\dots ,n_k}(2\pi i(n_1+\dots +n_k)\,z).
\eqno(3.13)
$$ 
Combining it with (3.12) and (3.3) we will get (3.8).

\smallskip

For $k=1$ we check (3.13) directly:
$$
(2\pi i)^{m_{v_1}}I^z_{i\infty}(\omega_{v_1})=
(2\pi i)^{m_{v_1}}\int_{i\infty}^z \omega_{v_1}(z_1)=
$$
$$
=(2\pi i)^{m_{v_1}}\sum_{n_1=1}^{\infty} c_{v_1,n_1} \int_{i\infty}^z
e^{2\pi in_1z_1} z_1^{m_{v_1}-1} dz_1.
$$
Putting in the $n_1$--th summand $t=2\pi in_1z_1$, we can rewrite this as
$$
\sum_{n_1=1}^{\infty} \frac{c_{v_1,n_1}}{n_1^{m_{v_1}}}\, 
\int_{-\infty}^{2\pi in_1z} e^tt^{m_{v_1}-1}dt\, .
\eqno(3.14)
$$
Since
$$
\int e^tP(t) dt =e^t (1+\partial_t)^{-1}P(t) + const \,,
\eqno(3.15)
$$
this is equivalent to (3.13).

\smallskip

The inductive step from $k$ to $k+1$ is similar: assuming (3.13)${}_k$, we have
$$
(2\pi i)^{m_{v_1}+\dots +m_{v_{k+1}}}I^z_{i\infty}(\omega_{v_{k+1}},\dots ,\omega_{v_1})=
$$
$$
(2\pi i)^{m_{v_{k+1}}}
\sum_{n_{k+1}=1}^{\infty} c_{v_{k+1},n_{k+1}} \int_{i\infty}^{z}
e^{2\pi in_{k+1}z_k} z_k^{m_{v_{k+1}}-1}  \times
$$
$$
\times \, \sum_{n_1,\dots ,n_k\ge 1}
\frac{c_{v_1,n_1}\dots c_{v_k,n_k}e^{2\pi i(n_1+\dots +n_k)z_k}}{n_1^{m_{v_1}}(n_1+n_2)^{m_{v_1}} \dots (n_1+\dots +n_k)^{m_{v_k}}}
\, D_{m_{v_1},\dots ,m_{v_k}}^{n_1,\dots ,n_k}(2\pi i(n_1+\dots +n_k)\,z_k)\,dz_k =
$$
$$
=(2\pi i)^{m_{v_{k+1}}}\sum_{n_1,\dots ,n_{k+1}\ge 1}
\frac{c_{v_1,n_1}\dots c_{v_k,n_k}c_{v_{k+1},n_{k+1}}}{n_1^{m_{v_1}}(n_1+n_2)^{m_{v_1}} \dots (n_1+\dots +n_k)^{m_{v_k}}}\,\times
$$
$$
\times \,\int_{i\infty}^z e^{2\pi i(n_1+\dots +n_k+n_{k+1})z_k} D_{m_{v_1},\dots ,m_{v_k}}^{n_1,\dots ,n_k}(2\pi i(n_1+\dots +n_k)\,z_k)\,z_k^{m_{v_{k+1}}-1} dz_k \,.
\eqno(3.16)
$$
Putting here $t=2\pi i(n_1+\dots +n_{k+1})\,z_k$, we can rewrite the last integral as
$$
\frac{1}{(2\pi i)^{m_{v_{k+1}}}}\,
\frac{1}{(n_1+\dots +n_{k+1})^{m_{v_{k+1}}}}\,
\int_{i\infty}^ze^t
D_{m_{v_1},\dots ,m_{v_k}}^{n_1,\dots ,n_k}
\left(\frac{n_1+\dots +n_k}{n_1+\dots +n_{k+1}}t\right)\,t^{m_{v_{k+1}}-1}dt 
$$
that is,
$$
\frac{1}{(2\pi i)^{m_{v_{k+1}}}}\,
\frac{1}{(n_1+\dots +n_k)^{m_{v_{k+1}}}}\,
e^{2\pi i(n_1+\dots +n_{k+1})z}\, \times
$$
$$
\times   \left.   \,(1+\partial_t)^{-1}(
D_{m_{v_1},\dots ,m_{v_k}}^{n_1,\dots ,n_k}
\left(\frac{n_1+\dots +n_k}{n_1+\dots +n_{k+1}}\,t\right)\,t^{m_{v_{k+1}}-1})\,
\right|_{t=2\pi i(n_1+\dots +n_{k+1})z}\,.
$$
Substituting this into (3.16), we finally obtain (3.13) and (3.12).

\medskip

{\bf 3.5. The limit $z\to 0.$} 
We can try to get an expression for
$$
I^0_{i\infty}(\omega_{v_k},\dots ,\omega_{v_1})
$$ 
as a (linear combination of)
multiple Dirichlet series, 
by formally putting $z=0$ in the r.h.s. of (3.8).
However, we will find out that this cannot be done automatically
for a certain range of values of $(j_k, j_{k-1})$, namely,
for $j_k\ge m_{v_k}+j_{k-1}-1-kC$: cf. Case 3 at the end of subsection
3.1.2.

\smallskip

To solve this problem, we will have for the first time to
assume that $z^{1-m_v}\omega_v(z)$ are of cusp modular type,
say, for the group $\Gamma_0(N)$, or for any modular subgroup 
which is normalized by the involution $z\mapsto gz$,
$$
g=g_N:=\left(\matrix 0& -1\\ N& 0\endmatrix \right) \,.
$$
We can then apply Proposition 2.4 which we reproduce and slightly augment:

\medskip

{\bf 3.5.1. Proposition.} {\it Assume that 
$\omega_V$ as above is a basis of a space of 1-forms
invariant with respect to $g_N$. Then
$$
J^0_{i\infty}(\omega_V)=
(g_{N*}(J^{\frac{i}{\sqrt{N}}}_{i\infty}(\omega_V)))^{-1} 
J^{\frac{i}{\sqrt{N}}}_{i\infty}(\omega_V)\,.
\eqno(3.17) 
$$
Replacing the coefficients of the formal series at the r.h.s of (3.17)
by their (convergent) representations via multiple Dirichlet series
(3.8), we get such representations for 
$I_{i\infty}^0(\omega_{v_k},\dots ,\omega_{v_1}) \,.$
}
\medskip

{\bf 3.5.2. Main application.} If we fix a modular subgroup
normalized by $g_N$, this proposition becomes applicable to
any family of cusp forms of the type $\omega_i(z) z^{m-1}$, $m\ge 1$, 
where $\omega_i(z)$ runs over a basis of the space of
forms of a fixed weight $2r$, and $m$ runs over $[1,2r-1]$ (cf. (2.4)).
Moreover, we can mix different weights, that is, take a finite
union of such families. 

\smallskip

Passing to a different basis of such a space, we may even assume that $(\omega_v)$
consists of eigenforms for $g_N$: $g^*_N(\omega_v )=
\varepsilon_v\omega_v$, $\varepsilon_v=\pm 1$, for all $s\in V.$ 

\smallskip

Coefficients of  $J^{\frac{i}{\sqrt{N}}}_{i\infty}$
are the series (3.3) at $z=\frac{i}{\sqrt{N}}$.

\bigskip

\centerline{\bf \S 4. Shuffle relations}

\smallskip

\centerline{\bf between multiple Dirichlet series}

\medskip

{\bf 4.1. Notation.} In this section, we will consider
(formal) multiple Dirichlet series of a special form
generalizing expressions (3.4), and deduce bilinear relations 
between them generalizing the well known {\it harmonic} shuffle relations
involving shuffles with repetitions.

\smallskip

Each such series will depend on a set of {\it coefficients data} $C$
and several complex or formal arguments $s_i$.
Here are precise definitions. Let $k\ge 1$ be a
natural number.

\medskip

{\bf 4.1.1. Definition.} {\it (i) Coefficients data $C$ of depth $k$
is a family of numbers $c_{n,m}^{(j,i)}$ indexed by two pairs of
integers satisfying $j>i\ge 0,\, j\le k$,
and  $n>m \ge 0$.

\smallskip

(ii) The multiple Dirichlet series associated with $C$
and arguments $s_1,\dots ,s_k$ is
$$
L_C(s_1,\dots , s_k):=
\sum_{0=u_0<u_1<\dots <u_k\in \bold{Z}}
\frac{\prod_{k\ge j>i \ge 0}c^{(j,i)}_{u_j,u_i}}{u_1^{s_1} u_2^{s_2} \dots u_k^{s_k}}
\eqno(4.1)
$$}
\smallskip

{\bf 4.1.2. Examples.} (a) Assume that $c_{n,m}^{(j,i)}=1$ if $m>0$ or $i>0$
and put  $c_{n,0}^{(j,0)}=a_n^{(j)}.$ Then
$$
L_C(s_1,\dots , s_k)= 
\sum_{0<u_1<\dots <u_k\in \bold{Z}}
\frac{a_{u_1}^{(1)}a_{u_1}^{(2)}\dots a_{u_k}^{(k)}}{u_1^{s_1} u_2^{s_2} \dots u_k^{s_k}}
\eqno(4.2)
$$
is an usual multiple Dirichlet series.

\medskip

(b) Define $c_{v,n}$ as in 3.1, and choose $v_1,\dots ,v_k\in V$
as in 3.1.1.  Construct the coefficients data $C$
putting 
$$
c^{(j,j-1)}_{n,m}:= c_{v_j,n-m},
\eqno(4.3)
$$
and  $c^{(j,i)}_{n,m}=1$ otherwise. Then 
$L_C(m_{v_1}+j_0-j_1,\dots , m_{v_k}+j_{k-1}-j_k)$
becomes the formal series (3.4) if we redenote $u_j=n_1+\dots +n_j.$

\medskip

{\bf 4.1.3. Shuffles and a composition of the coefficients data.}
Let $C=(c^{(j,i)}_{n,m})$ and $D=(d^{(j,i)}_{n,m})$ be two data
of depths $k$ and $l$ respectively. {\it A $(k,l,p)$--shuffle 
with repetitions}
is a pair of strictly increasing maps $\sigma =(\sigma_1,\sigma_2)$,
$$
\sigma_1:\, [0,k]\to [0,p],\quad \sigma_2:\, [0,l]\to [0,p]
$$
satisfying the following conditions:
$$
\sigma_1(0)=\sigma_2(0)=0, \ \sigma_1 ([0,k])\cup \sigma_2 ([0,l]) =[0,p] .
\eqno(4.4)
$$
It follows that $\roman{max}\,(k,l)\le p \le k+l.$ We will say
that {\it the $\sigma$--multiplicity of $j\in [0,p]$ is one},
if $j\notin \sigma_1 ([0,k])\cap \sigma_2 ([0,l])$.
Otherwise  the $\sigma$--multiplicity of $j\in [0,p]$ is two. 
In particular,  the $\sigma$--multiplicity of $0$ is two.

\smallskip

Given such $C$, $D$ and $\sigma =(\sigma_1,\sigma_2)$, we will define
the third coefficients data $E=(e^{(j,i)}_{n,m})$ of depth $p$
which we denote $E=C*_{\sigma}D$.  Choose $j,i$ with 
$p\ge j>i \ge 0.$ We have the following
set of mutually exclusive and exhausting alternatives
(A)${}_1$, (A)${}_2$, (B), (C).

\smallskip

(A) Assume that both $j$ and $i$ have multiplicities one.
\smallskip

(A)${}_1$. Both $j$ and $i$ belong to the image
of one and  the same $\sigma_a$ with $a=1$ or $a=2$. Then we put
$$
e^{(j,i)}_{n,m}:=c^{(\sigma_1^{-1}(j),\sigma_1^{-1}(i))}_{n,m} \
\roman{for}\ a=1,
$$
and
$$
e^{(j,i)}_{n,m}:=d^{(\sigma_2^{-1}(j),\sigma_2^{-1}(i))}_{n,m} \
\roman{for}\ a=2.
$$

(A)${}_2$. Assume that $j$, resp. $i$, belongs to the image
of $\sigma_a$, resp. $\sigma_b$, with $a\ne b.$ 
Then we put
$$
e^{(j,i)}_{n,m}:=1.
$$

\smallskip

(B) Assume that exactly one of $j,i$ has multiplicity two.
Then there exists only one value $a=1$ or $2$ such that
$j$ and $i$ belong to the image
of $\sigma_a$. We put then as in the case (A)${}_1$  
$$
e^{(j,i)}_{n,m}:=c^{(\sigma_1^{-1}(j),\sigma_1^{-1}(i))}_{n,m} \
\roman{for}\ a=1,
$$
and
$$
e^{(j,i)}_{n,m}:=d^{(\sigma_2^{-1}(j),\sigma_2^{-1}(i))}_{n,m} \
\roman{for}\ a=2.
$$
\smallskip

(C) Assume that both $i$ and $j$ have multiplicities two. Then
we put
$$
e^{(j,i)}_{n,m}:=c^{(\sigma_1^{-1}(j),\sigma_1^{-1}(i))}_{n,m}
d^{(\sigma_2^{-1}(j),\sigma_2^{-1}(i))}_{n,m} .
$$ 

\medskip

{\bf 4.1.4. Shuffles and a composition of the arguments.}
Let $s:=(s_1,\dots ,s_k)$ and $t:=(t_1,\dots ,t_l)$
be arguments for the data $C$ and $D$ as above,
and $\sigma$ a $(k,l,p)$--shuffle as above.
We define $s+_{\sigma}t :=(r_1,\dots ,r_p)$ as follows.

\smallskip

If $i$ has multiplicity one and is covered by $\sigma_1$, resp.
$\sigma_2$,
then $r_i:=s_{\sigma_1^{-1}(i)}$, resp.
$r_i:=t_{\sigma_2^{-1}(i)}$.

\smallskip

If $i$ has multiplicity two, then $r_i:=s_{\sigma_1^{-1}(i)}+
t_{\sigma_2^{-1}(i)}$.

\smallskip

We can now state the main result of this section.

\medskip

{\bf 4.2. Theorem.} {\it Let $C$, resp. $D$,
be some coefficients data of depths $k$, resp. $l$, as above.
Then we have
$$
L_C(s)\cdot L_D(t)=\sum_{\sigma} L_{C*_{\sigma}D}(s+_{\sigma}t)
\eqno(4.5)
$$
where the summation is taken over all $(k,l,p)$--shuffles with repetitions.}

\smallskip

{\bf Proof.} Consider a term of the series (4.1)
corresponding to $(u_0=0,u_1,u_2,\dots ,u_k)$
and a term of the series $L_D(t)$ corresponding
to, say, $(w_0=0, w_1, w_2,\dots ,w_l).$ This pair of terms
determines a unique $(k,l,p)$--shuffle $(\sigma_1,\sigma_2)$, where
$p$ is the cardinality of the union of sets 
$$
\{u_1,u_2,\dots ,u_k\}\cup
\{ w_1, w_2,\dots ,w_l\}:= \{q_1, \dots ,q_p\}.
$$ 
Namely, we may and will assume that $q_0=0<q_1<\dots <q_p$.
Then $\sigma_1(i)=j$ if $u_i=q_j$, and $\sigma_2(i)=j$ if $w_i=q_j$.

\smallskip  

Group together all pairwise products corresponding to one and the same
shuffle, and denote the resulting sum by $L_{\sigma}$.

\smallskip

The denominator of one such a product will obviously be
$q_1^{r_1}\dots q_p^{r_p}$ where $r=s+_{\sigma}t$. 
Moreover, knowing such a denominator, we uniquely reconstruct
the two terms from $L_C(s)$ and $L_D(t)$, from which it was
produced, at least if $s,t$ take generic values so that
in the family $\{ s_a, t_b, s_a+t_b\}$ all terms are
pairwise distinct.
Finally, all possible sequences $q_0=0<q_1<\dots <q_p$ 
will occur.

\smallskip

To prove that $L_{\sigma}=L_{C*_{\sigma}D}(s+_{\sigma}t)$, it remains to check that
the numerator of such a product will be as predicted by (4.5),
in other words, that
$$
\prod _{p\ge j>i\ge 0}e_{q_j,q_i}^{(j,i)} = (?)
\prod _{k\ge j>i\ge 0}c_{u_j,u_i}^{(j,i)}
\prod _{l\ge j>i\ge 0}d_{w_j,w_i}^{(j,i)}
$$
if $e^{(j,i)}_{n,m}$ are defined as in 4.1.3.

\smallskip

This is straightforward, although somewhat tedious.

\medskip

{\bf 4.3. Concluding remarks.} It would be interesting
to describe some nontrivial spaces of Dirichlet series
containing periods of cusp forms, closed with respect
to the series shuffle relations, and consisting entirely
of periods in the sense of [KonZa].

\smallskip

Regarding shuffle relations themselves, motivic philosophy
predicts that they should be obtainable by standard manipulations
with integrals. For the harmonic shuffle relations between
multiple zeta values,
A.~Goncharov established this in the Ch. 9 of
[Go6] (for convergent integrals), and in 7.5
of [Go6] elaborating the last page of [Go4] (for regularized
integrals). Conversely, integral shuffle relations
can be deduced from harmonic ones: see [Go5], Ch. 2.

\newpage

\centerline{\bf \S 5. Iterated Eichler--Shimura and Hecke  relations}

\medskip

{\bf 5.1. Eichler--Shimura relations for iterated integrals.} In this subsection, we take for $X$ the upper
half--plane $H$ and for $(\omega_V)$ a family of 1--forms 
of cusp modular type (see 2.1.1) spanning
a finite--dimensional linear space stable with respect to 
the modular transformations
$$
\sigma = \left(\matrix 0& -1\\1 &0\endmatrix \right) \,, \quad
\tau  = \left(\matrix 1& -1\\1 &0\endmatrix \right) \,.
$$

\smallskip

{\bf 5.1.1. Proposition.} {\it With these assumptions, we have
$$
 \sigma_* (J_0^{i\infty}(\omega_V))\cdot J_0^{i\infty}(\omega_V) =1.
\eqno(5.1)
$$
$$
\tau_*^2(J_0^{i\infty}(\omega_V))\cdot \tau_*(J_0^{i\infty}(\omega_V))\cdot
J_0^{i\infty}(\omega_V)=1.
\eqno(5.2)
$$
}  
{\bf Proof.} From (1.15) we get
$$
\sigma_* (J_0^{i\infty}(\omega_V))=J_{\sigma (0)}^{\sigma (i\infty)}(\omega_V)=
J^0_{i\infty}(\omega_V)=J_0^{i\infty}(\omega_V)^{-1}.
$$
which shows (5.1).

\smallskip

Similarly, $\tau$ transforms $(0,i\infty )$ into
$(i\infty , 1)$, then to $(1,0)$, and according to (1.9)
$$
J_1^0(\omega_V)\cdot J^1_{i\infty}(\omega_V)\cdot
J_0^{i\infty}(\omega_V)
=1
\eqno (5.3)
$$
which is the same as (5.2).

\smallskip

Notice that some care is needed in establishing (5.3):
the geodesic triangle with vertices $0, i\infty, 1$ should be first
replaced by a sequence of geodesic hexagons lying entirely
in $H$ and cutting the corners of  the triangle,
and then it must be checked that in the limit
the hexagon relation replacing (5.3) tends to (5.3).
This is routine for cusp modular 1--forms, cf. 2.1.3.

\smallskip

We now pass to the relations involving Hecke operators.

\medskip

{\bf 5.2. Hecke operators.} In this subsection,
$V$ denotes  a finite set, $\omega_v = f_v(z)dz$
a family of modular forms {\it of weight two},
and $p_v$ a family of primes, both indexed by $v\in V$.
Moreover, we assume that $T_{p_v}\omega_v =\lambda_v\omega_v$,
$\lambda_v\in \bold{C}$,
where $T_{p_v}$ is the Hecke operator 
$$
T_{p_v}:= \left(\matrix p_v& 0\\ 0& 1\endmatrix \right) +
\sum_{b=0}^{p_v-1}
\left(\matrix 1& b\\ 0& p_v\endmatrix \right) =\bold{p_v}+
\sum_{b=0}^{p_v-1}h(p_v,b)\,.
\eqno(5.4)
$$

\smallskip

Put $U:=V\coprod V^{\prime}$ where $V^{\prime}$ is another copy of 
the indexing set $V$, and for $v^{\prime}\in  V^{\prime}$
corresponding to $v\in V$ put  $\omega_{v^{\prime}}:=(\bold{p_v})^*
(\omega_v)$. Let $\omega_U$ be the family consisting of
all $\omega_v$ and $\omega_{v^{\prime}}$. When we consider
formal series of the type $J_a^z(\omega_U)$ as in \S 1,
we denote the variables corresponding to $V$, resp. $V^{\prime}$,
by $A_v$, resp. $A_{v^{\prime}}.$

\smallskip

Denote by $W$ the set of pairs $w=(v, b)$ where $v\in V$ and
$b\in [0,p_v-1]$. Let $\omega_W$ be the family
consisting of $\omega_{(v, b)}:=h(p_v,b)^*(\omega_v)$.
When we consider
formal series of the type $J_a^z(\omega_W)$,
we denote the variables corresponding to $w$,
by $B_w$ or $B_{(v,b)}$.

\smallskip

Define the following two continuous homomorphisms of
rings of formal series:
$$
l:\, \bold{C}\langle\langle A_U\rangle\rangle \to
\bold{C}\langle\langle A_V\rangle\rangle\,:\quad
l(A_v):=\lambda_vA_v,\ l(A_{v^{\prime}}):=-A_v\,,
\eqno(5.5)
$$
$$
r:\, \bold{C}\langle\langle B_W\rangle\rangle \to
\bold{C}\langle\langle A_V\rangle\rangle\,:\quad
r(B_{(v,b)}):=A_v \,.
\eqno(5.6)
$$

\medskip

{\bf 5.3. Theorem.} {\it We have
$$
l(J_{i\infty}^0(\omega_U))=r(J_{i\infty}^0(\omega_W))\,.
\eqno(5.7)
$$
}
\smallskip

{\bf Proof.} We will check that 
$$
l(J_{i\infty}^0(\omega_U)) =J^0_{i\infty}((\lambda_v\omega_v-
\omega_{v^{\prime}}))
\eqno(5.8)
$$
whereas
$$
r(J_{i\infty}^0(\omega_W)) =J^0_{i\infty}((\sum_{b=0}^{p_v-1}\omega_{(v,b)}))\,
\eqno(5.9)
$$
where at the right hand sides we consider both families as indexed by $V$.
Since from (5.4) and the definitions above 
we obtain for each $v\in V$
$$
\lambda_v\omega_v-\omega_{v^{\prime}}=\sum_{b=0}^{p_v-1}\omega_{(v,b)}.
\eqno(5.10)
$$
this will prove the theorem.

\smallskip 

We have
$$
J_{i\infty}^0(\omega_U)=
\sum_{n=0}^{\infty}\sum_{(u_1,\dots ,u_n)\in U^n}
A_{u_n}\dots A_{u_1} I_{i\infty}^0(\omega_{u_n},\dots ,\omega_{u_1})\, .
\eqno(5.11)
$$
Consider one summand in (5.11). In the sequence 
$(u_1,\dots ,u_n)$ there are several, say $0\le k\le n$,
elements $v_i\in V$ and the remaining $n-k$ elements
$v^{\prime}_j\in V^{\prime}$. Application of
$l$ eliminates all primes in the subscripts
of $A_{u_n}\dots A_{u_1}$ and produces a monomial
in $A_v$; besides, it multiplies this monomial
by $(-1)^{n-k}\prod \lambda_{v_i}$. Hence the coefficient
at any monomial $A_{v_n}\dots A_{v_1}$ in $l(J_{i\infty}^0(\omega_U))$
can be written as
$$
\sum_{S\subset [1,\dots n]} (-1)^{n-|S|}\left(\prod_{i\in S} \lambda_{v_i}
\right) I_{i\infty}^0(\pi_S (\omega_{v_n},\dots ,\omega_{v_1}))
\eqno(5.12)
$$
where the operator $\pi_S$ replaces $v_j$ by $v^{\prime}_j$
whenever $j\notin S.$

\smallskip
On the other hand, the similar coefficient in
$J^0_{i\infty}((\lambda_v\omega_v-
\omega_{v^{\prime}}))$ is
$$
I_{i\infty}^0( \lambda_{v_n}\omega_{v_n}-\omega_{v_n^{\prime}},
\dots , \lambda_{v_1}\omega_{v_1}-\omega_{v_1^{\prime}})
$$
which obviously coincides with (5.12) because the iterated integrals
are polylinear in $\omega$. This proves (5.8).

\smallskip

The check of (5.9) is similar. We have
$$
J^0_{i\infty}(\omega_W)=
$$
$$
=\sum_{((v_n,b_n),\dots ,(v_1,b_1))\in W^n}
B_{(v_n,b_n)} \dots B_{(v_1,b_1)}
 I_{i\infty}^0(\omega_{(v_n,b_n)},\dots ,
\omega_{(u_1,b_1)})\, .
$$
Application of $r$ produces a series in $(A_v)$
whose coefficient at   $A_{v_n}\dots A_{v_1}$ equals
$$
\sum_{(b_1,\dots ,b_n)} I_{i\infty}^0(\omega_{(v_n,b_n)},\dots ,
\omega_{(v_1,b_1)})\, .
$$
This is the same as the respective coefficient at the r.h.s. of (5.9).

\bigskip

\centerline{\bf \S 6. Differentials of the third kind}

\smallskip

\centerline{\bf and generalized associators}

\medskip

{\bf 6.1. Normalized horizontal sections.} In this
subsection, following [Dr2], we will define and study solutions of the differential
equation (1.4), $dJ^z(\Omega )=\Omega (z)\,J^z(\Omega )$
in the case when $\Omega = \sum A_v\omega_v$ may
have a logarithmic singularity at  a point $a$ so that
it cannot be normalized by the condition $J^a(\Omega )=1$
and in fact cannot be defined by the series (1.3).
   
\smallskip

We start with a local situation. Put 
$$
r_{v,a} :=\roman{res}_a\, \omega_v,\quad R_a:= \roman{res}_a\,\Omega =
\sum_v r_{v,a} A_v\,.
\eqno(6.1)
$$
The normalized solution will depend on the choice of a local parameter
$t_a$ at $a$, and a branch of logarithm $\roman{log}\,t_a$.

\smallskip
 
Let $U$ be a disc around $a$
uniformized by $t_a$. Denote by $\roman{log}\,t_a$ the branch
of the logarithm in $U$
which is real on $\roman{Im}\,t_a=0,\, \roman{Re}\,t_a>0.$
Delete from $U$ a cut from $a$ to the boundary which does not
intersect the latter interval and denote $U^{\prime}$
the remaining domain.
Write $t_a^{R_a}$ for $e^{R_a \roman{log}\,t_a}$. It is a formal
series in $A_v$ with coefficients which are holomorphic
functions near $a$ in $U^{\prime}$. Assume that outside of $a$, all $\omega_v$ 
are regular in $U$. 

\medskip

{\bf 6.1.1. Definition.} {\it A $\nabla_{\Omega}$--horizontal section
$J$  in $U^{\prime}$ is 
called normalized at $a$ (with respect to a choice of $t_a$)
if it is of the form $J=K\cdot t_a^{R_a}$, where $K$ can be extended
to a holomorphic section in some neighborhood of $a$ in $U$
which takes the value $1$ at $a$.}

\smallskip

We will see that this definition produces a version of $J_a^z(\Omega )$.
In fact,  we get precisely $J_a^z(\Omega )$, if $R_a=0$ so that $t_a^{R_a}=1$.

\medskip

{\bf 6.1.2. Proposition.} {\it For any $a$ and $t_a$ as above, there exists 
a unique normalized at $a$ local section holomorphic in $U^{\prime}$.}

\smallskip

{\bf Proof.} In the course of proof, we will be considering
only one point $a$, so we will omit it in the notation for brevity
and write $R,t, r_v$ etc in place of former $R_a, t_a, r_{v,a}$.

\smallskip

The equation $\nabla_{\Omega} (K\cdot t^R)=0$ is equivalent to
$$
dK=\Omega^{\prime}K + t^{-1}[R,K]\,dt\,
\eqno(6.2)
$$
where 
$$
\Omega^{\prime}:=\Omega - R\dfrac{dt}{t} =\sum_v r_v\nu_v A_v\,.
$$
We look for a solution to (6.2) of the form 
$$
K= 1+\sum_{n=1}^{\infty}\sum_{(v_1,\dots ,v_n)\in V^n} 
f_{v_1,\dots ,v_n} A_{v_1}\dots A_{v_n}
$$ 
where $f_{v_1,\dots ,v_n}$ must be holomorphic functions
defined in some common neighborhood of $a$ in $U$
and vanishing at $a$. From (6.2) we find $df_v=\nu_v$
so that the only choice is $f_v(z):=\int_a^z \nu_v$. Notice that all $\nu_v$ are regular in a common
neighborhood of $a$ uniformized by $t_a$. We will see that in this neighborhood
all other coefficients can be defined as well. 

\smallskip

If $f_{v_1,\dots ,v_{n-1}}$ with required properties are defined for
some $n-1\ge 1$, we find from (6.2)
$$
df_{v_1,\dots ,v_n}=\nu_{v_1} f_{v_2,\dots ,v_n}
+t^{-1}(r_{v_1} f_{v_2,\dots ,v_n} -  f_{v_2,\dots ,v_{n-1}}r_{v_n})\,dt.
$$
By the inductive assumption, the r.h.s. is well defined and regular at $a$,
so integrating it from $a$ to $z$ we get $f_{v_1,\dots ,v_n}(z)$.

\medskip

{\bf 6.2. Scattering operators.}  Let now $a,b$ be two
points where $\Omega$ may have
logarithmic singularities. Choose $t_a$, $t_b$ as above,
construct the neighborhoods $U_a$, $U_b$ and neighborhoods with deleted
cuts $U_a^{\prime}$, $U_b^{\prime}$ in which we have the holomorphic
normalized horizontal sections $J_a$, $J_b$. Now
embed $U_a^{\prime}$, $U_b^{\prime}$ into a connected simply connected domain
$W$ to which both $J_a$ and $J_b$ can be analytically extended.
Clearly, they are invertible elements of $\Cal{O}_x\langle\langle A_V
\rangle\rangle$ at almost all points $x\in W$.
Put
$$
\widetilde{J}_b^a = J_a^{-1}J_b \,.
\eqno(6.3)
$$
As in the proof of Proposition 1.2, one sees that 
$\widetilde{J}_a^b  \in \bold{C}\langle\langle A_V\rangle\rangle$.
Borrowing the physics terminology, we can call this transition
element {\it the scattering operator}. 

\smallskip

I added twiddle in the notation in order to remind
the reader that $\widetilde{J}_a^b$, besides $\Omega$, depends on $t_a$
and $t_b$ as well, if at least one of the residues $R_a, R_b$
does not vanish. This dependence however is
pretty mild. Let  $J_a^{\prime}=K^{\prime} (t_a^{\prime})^{R_a}$ be another horizontal section
normalized with respect to some $t_a^{\prime}.$ Denote by
$\tau_a \in \bold{C}^*$ the value of $t_a^{\prime}/t_a$ at $a$,
and let $t_a^{\prime} =  T_a\cdot t_a\tau_a$, $T_a(a)=1$. 

\medskip

{\bf 6.2.1. Proposition.} {\it (i) We have
$$
J_a^{\prime} = J_a \cdot \tau_a^{R_a}, \quad K=K^{\prime}\cdot T^{R_a}\,.
\eqno(6.4) 
$$
(ii) Therefore, after replacing two uniformizing
parameters $t_a, t_b$  by $t_a^{\prime}, t_b^{\prime}$,
we get
$$
\widetilde{J}_b^{\prime a} = \tau_a^{-R_a}\widetilde{J}_b^a 
\tau_b^{R_b}\,.
\eqno(6.5)
$$
} 
\smallskip

{\bf Proof.} (i) We have
$$
J_a^{\prime}=K^{\prime} (t_a^{\prime})^{R_a}=
K^{\prime}\cdot T_a^{R_a}\cdot  (t_a)^{R_a}\tau^{R_a}
$$
Since $\tau^{R_a}  \in \bold{C}\langle\langle A_V\rangle\rangle$,
and this section is invertible,
$K^{\prime}\cdot T_a^{R_a}\cdot  (t_a)^{R_a}$ is
$\nabla_{\Omega}$--horizontal as well,
and since $K^{\prime}\cdot T_a^{R_a}$ at $a$ equals 1,
it is normalized, so  $K^{\prime}\cdot T_a^{R_a}\cdot  (t_a)^{R_a}=J_a$, which proves (6.4).

\smallskip

(ii) Similarly, from $J_b^{\prime} = J_a^{\prime} \widetilde{J}_b^{\prime a}$
and (6.4) we get
$$
J_b=J_a\,\tau_a^{R_a} \, \widetilde{J}_b^{\prime a}\, \tau_b^{-R_b}
$$
which together with (6.3) proves (6.5).

\medskip

{\bf 6.3. Example: Drinfeld's associator.} Let $X=\bold{P}^1(\bold{C})$,
$V=\{0,1\}$,
$$
\omega_0=\frac{1}{2\pi i}\,\frac{dz}{z},\quad
\omega_1=\frac{1}{2\pi i}\,\frac{dz}{z-1}\,.
$$
Then
$$
\Omega = A_0\omega_0 + A_1\omega_1
$$
has poles at $0,1,\infty$ with residues $A_0/2\pi i$, $A_1/2\pi i$,
$-(A_0 + A_1)/2\pi i$ respectively. Put $t_0=z, t_1= 1-z$.
Then $\widetilde{J}_0^1$ in our notation is
the Drinfeld associator  $\phi_{KZ}(A_0,A_1)$ from \S 2 of [Dr2].

\medskip

{\bf 6.4. Generalized associators.} An essential feature of the
the last example is that $\omega_v$ are global logarithmic
differentials on the compact Riemannian surface $\bold{P}^1$.

\smallskip

Generally, let $X$ be such a surface, and $(a_i)$ a finite set of 
$N$ points on it. The dimension of the space of global
logarithmic differentials with poles in this set, $\sum c_i\, d\,\roman{log}\,f_i$, where
$c_i\in\bold{C}$, $f_i$ are meromorphic on $X$, 
is bounded by $N-1$. It achieves the maximum value
$N-1$ iff the difference of any two points $a_i-a_j$
is torsion in the divisor class group.  I will call
such a set $(a_i)$  {\it logarithmic.} The supply of logarithmic set
depends on the genus of $X$.

\smallskip

a) {\it Genus zero.} Any finite set of points is logarithmic.
The respective iterated integrals in this case include
multiple polylogarithms introduced in [Go1]. In fact, as Goncharov
remarked, general iterated integrals can be reduced
to multiple polylogarithms.

\smallskip

b) {\it Genus one.} In this case we can take any finite set of points 
of finite order, for example, the subgroup of all points of a given order
$M$.

\smallskip

For a general subset of points, Goncharov found a Feynman integral
presentation (in the sense of the last section of [Go2])
of the respective real periods. His formulas involve
the generalized Kronecker--Eisenstein series and can be considered
as an elliptic version of multiple zeta values. A similar
problem is addressed in the recent work of A.~Levin and G.~Racinet.

\smallskip
 
c) {\it Higher genera.} If genus of $X$ is $>1$, the order
of such a logarithmic set is bounded. The most interesting explicitly known
examples are modular curves and cusps on them, cf. [Dr1], [El].

\smallskip

Notice, that the initial Drinfeld's setting is modular
as well: $\bold{P}^1$ with three marked points ``is''
the modular curve $\Gamma_0(4)\setminus \overline{H}$
together with its cusps.

\smallskip

According to Deligne and Elkik ([El]), a set of points is logarithmic iff
the mixed Hodge structure on $H^1(X^{\circ},\bold{Q})$ (where $X^{\circ}$
is the complement to ) is split that is, the direct sum of pure Hodge structures.

\medskip

{\bf 6.4.1. Definition.} {\it Let $X$ be a compact
Riemannian surface, and $(a_i)$ a logarithmic set of points on it.
Then any scattering operator of the form $\widetilde{J}^{a_j}_{a_i}$
is called a generalized associator.}

\medskip

{\bf 6.5. Relations between scattering operators.} Three types
of relations established for $J_a^z(\omega_V)$ in \S 1
are extended below to the case of the general scattering operators.

\medskip
 
{\bf 6.5.1. Group--like property.} {\it We have
$$
\Delta (\widetilde{J}_a^b) = \widetilde{J}_a^b
\otimes \widetilde{J}_a^b 
\eqno(6.6)
$$
where $\Delta$ is defined in 1.4.1.}

\smallskip

To see this, notice that $\Delta (J_a)$, by definition,
is the series which is obtained from $J_a$
by replacing each $A_v$ with $B_v:=A_v\otimes 1 +1\otimes A_v$.
Hence $\Delta (R_a) = R_a\otimes 1 + 1\otimes R_a$ is the residue
of $\Delta (\Omega )$ at $a$. Therefore
$\Delta (J_a)$ is the normalized $\nabla_{\Delta(\Omega)}$--horizontal section 
in the ring of formal series in $B_v$, and $\Delta (\widetilde{J}_a^b)$
is the respective scattering operator. 

\smallskip

On the other hand, $J_a\otimes J_a$ satisfies the same equation
$d\,(J_a\otimes J_a) = (\Omega\otimes 1+1\otimes\Omega )(J_a\otimes J_a)$
and clearly is as well normalized at $t_a$. Hence the passage
from $J_a\otimes J_a$ to $J_b\otimes J_b$ is governed by the same
scattering operator, this time represented as the r.h.s. of (6.6).

\medskip

{\bf 6.5.2. Cycle identities.} {\it Let $\gamma$ be a closed oriented contractible contour in $U$, inside which there are no singularities of $\Omega$.
Let $a_1, \dots ,a_n$ be points along this contour
(cyclically) ordered compatibly with orientation. Then
$$
\widetilde{J}^{a_1}_{a_2}\widetilde{J}^{a_2}_{a_3}\dots 
\widetilde{J}^{a_{n-1}}_{a_n}
\widetilde{J}^{a_n}_{a_1} =1\,.
\eqno(6.7)
$$
}

Of course, in this statement we assume that at each point 
$a_i$ one and the same 
local parameter $t_{a_i}$ and one and the same branch
of its logarithm is used for the  definition of the two relevant
normalized sections corresponding to the incoming and outcoming
segments of the contour. Otherwise the relevant $\tau^R$
factors as in (6.5) must be inserted.

\medskip

{\bf 6.5.3. Functoriality.} {\it Let $g$ be an automorphism
of $X$ compatibly acting upon all the relevant objects, with
the possible exception of 
parameters $t_a$, and transforming the space
spanned by $\omega_v$ into itself. Define $g_*$ as in 1.4. Then
$$
\widetilde{J}_{ga}^{gb} =\tau_a^{-R_a}g_*(\widetilde{J}_a^b)\, \tau_b^{R_b} \,.
 \eqno(6.8)
$$
where the $\tau^R$ factors account for the passage
from $t_a$, resp. $t_b$, to $g^*(t_a)$, resp. $g^*(t_a)$.   
}

The proof is essentially the same as in 1.4 and 1.5.

\medskip

{\bf 6.6. Example: Drinfeld associator revisited.}
If we treat Drinfeld's setup as the $\Gamma_0(4)$--modular
curve, lift it to $H$ and apply to the respective family
of scattering operators the Eichler--Shimura
relations (5.1), (5.2) (following from the cycle identities
of lengths 2 and 3, and
functoriality), we will get the duality and
the hexagonal relations. Of course, this is how they were deduced
in the first place, with $\sigma ,\tau$ pushed down to 
$\bold{P}^1$ rather than everything else lifted to $H$.

\smallskip

It seems very probable that the somewhat mysterious
relationship between the double logarithms at roots of unity
and the modular complex discovered in [Go3]
can be explained in the same way.

\bigskip

\centerline{\bf References}

\medskip

[AB] A.~Ash, A.~Borel. {\it Generalized modular symbols.}
In: Lecture Notes in Math., 1447 (1990), Springer Velag, 57--75.

\smallskip

[AR] A.~Ash, L.~Rudolph. {\it The modular symbol and
continued fractions in higher dimensions.} Inv. Math. 55 (1979), 241--250.

\smallskip

[Ch] K.--T. Chen. {\it Iterated path integrals.}
Bull. AMS, 83 (1977), 831--879.

\smallskip

[De] P.~Deligne. {\it Multizeta values.} Notes d'expos\'es, IAS, Princeton, 2001.

\smallskip

[DeGo]  P.~Deligne A.~Goncharov. {\it Groupes fondamentaux
motiviques de Tate mixte.} e--Print math.NT/0302267 .

\smallskip

[Dr1] V.~G.~Drinfeld. {\it Two theorems on modular curves.} Func. An.
and its Applications, 7:2 (1973), 155--156.

\smallskip

[Dr2] V.~G.~Drinfeld. {\it On quasi--triangular quasi--Hopf
algebras and some groups closely associated
with $Gal\,(\overline{\bold{Q}}/\bold{Q})$.}
Algebra and Analysis 2:4 (1990); Leningrad Math. J.
2:4 (1991), 829--860. 

\smallskip

[El] R. Elkik. {\it Le th\'eor\`eme de Manin--Drinfeld.}
Ast\'erisque 183 (1990), 59--67.

\smallskip

[Go1] A.~Goncharov. {\it Polylogarithms in arithmetic and geometry.}
Proc. of the ICM 1994 in Z\"urich,  vol. 1, Birkh\"auser, 1995, 374--387.

\smallskip

[Go2] A.~Goncharov. {\it Multiple $\zeta$-values, Galois
groups and geometry of
modular varieties}. Proc. of the Third European Congress of
Mathematicians. Progr. in Math.,
vol. 201, p. 361--392. Birkh\"auser Verlag, 2001.
e--Print math.AG/0005069.

\smallskip

[Go3] A.~Goncharov. {\it The double logarithm and Manin's complex for modular curves.} Math. Res. Letters, 4 (1997), 617--636.

\smallskip

[Go4] A.~Goncharov. {\it Multiple polylogarithms,
cyclotomy, and modular complexes.}  Math. Res. Letters, 5 (1998), 497--516.

\smallskip

[Go5] A.~Goncharov. {\it Multiple polylogarithms and mixed Tate motives.}
e--Print math.AG/0103059

\smallskip

[Go6] A.~Goncharov. {\it Periods and mixed motives.}
e--Print math.AG/0202154

\smallskip

[GoMa] A.~Goncharov, Yu.~Manin. {\it Multiple zeta-motives and moduli spaces $\overline{M}_{0,n}$.} Compos. Math.
140:1 (2004), 1--14.  e--Print math.AG/0204102

\smallskip

[Ha] R.~Hain. {\it Iterated integrals and algebraic cycles:
examples and prospects.} e--Print math.AG/0109204

\smallskip

[Kon] M.~Kontsevich. {\it Operads and motives in deformation
quantization.} Lett. Math. Phys. 48 (1999), 35--72.

\smallskip

[KonZa] M.~Kontsevich, D.~Zagier. {\it Periods.} In: Mathematics
Unlimited --- 2001 and Beyond, Springer, Berlin, 2001, 771--808.

\smallskip

[Ma1] Yu.~Manin. {\it Parabolic points and zeta-functions of modular curves.}
 Russian: Izv. AN SSSR, ser. mat. 36:1 (1972), 19--66. English:
Math. USSR Izvestija, publ. by AMS, vol. 6, No. 1 (1972), 19--64,
and Selected papers, World Scientific, 1996, 202--247.

\smallskip

[Ma2] Yu.~Manin. {\it Periods of parabolic forms and $p$--adic Hecke series.}
Russian: Mat. Sbornik, 92:3 (1973), 378--401. English:
Math. USSR Sbornik, 21:3 (1973), 371--393
and Selected papers, World Scientific, 1996, 268--290.

\smallskip

[MaMar] Yu.~Manin, M.~Marcolli. {\it Continued fractions, modular symbols, and non-commutative geometry.} Selecta math., new ser. 8 (2002),
475--521. e--Print math.NT/0102006

\smallskip

[Me] L.~Merel. {\it Quelques aspects arithm\'etiques et g\'eom\'etriques
de la th\'eorie des symboles modulaires.} Th\'ese de doctorat,
Universit\'e Paris VI, 1993.

\smallskip

[Ra1] G.~Racinet. {\it S\'eries g\'en\'eratrices non--commutatives
de polyz\^etas et associateurs de Drinfeld.}
Th\`ese de doctorat, http://www.dma.ens.fr/$\sim$racinet

\smallskip

[Ra2] G.~Racinet. {\it Doubles m\'elanges des polylogarithms multiples aux racines de l'unit\'e.} Publ. Math. IHES, 95 (2002), 185--231.
e--Print math.QA/0202142

\smallskip

[Ra3] G.~Racinet. {\it Summary of algebraic relations between
multiple zeta values.} Notes of the talk at MPIM, Bonn,
Aug. 17, 2004.

\smallskip

[Sh1] V.~Shokurov. {\it Modular symbols of arbitrary weight.} Func. An.
and its Applications, 10:1 (1976), 85--86.

\smallskip

[Sh2] V.~Shokurov. {\it The study of the homology of Kuga
varieties.} Math. USSR Izvestiya, 16:2 (1981), 399--418.

\smallskip

[Sh3] V.~Shokurov. {\it Shimura integrals of cusp forms.}
Math. USSR Izvestiya, 16:3 (1981), 603--646.

\smallskip

[T] T.~Terasoma. {\it Mixed Tate motives and multiple zeta values.}
Inv. Math. 149 (2002), 339--369.

\smallskip

[Za1] D.~Zagier. {\it Hecke operators and periods of modular forms.}
In: Israel Math. Conf. Proc., vol. 3 (1990), I.~I.~Piatetski--Shapiro
Festschrift, Part II, 321--336.

\smallskip

[Za2] D.~Zagier. {\it Values of zeta functions and their
applications.} First European Congress of Mathematics,
vol. II, Birkh\"auser Verlag, Basel, 1994, 497--512.

\smallskip

[Za3] D.~Zagier. {\it Periods of modular forms and Jacobi
theta functions.} Inv. Math., 104 (1991), 449--465.

\enddocument